\documentclass[a4paper,10pt]{article}

\usepackage{amssymb}
\usepackage{amsmath}
\usepackage{amsthm}
\usepackage{amscd}
\usepackage[english]{babel}
\usepackage[all]{xypic}
\usepackage[all]{xy}
\usepackage{url}
\usepackage{graphicx}
\usepackage{hyperref}

\usepackage{comment}

\newtheorem{theo}{Theorem}[section]
\newtheorem*{theo*}{Theorem}
\newtheorem{prop}[theo]{Proposition}
\newtheorem{lem}[theo]{Lemma}
\newtheorem{cor}[theo]{Corollary}
\newtheorem{defi}[theo]{Definition}
\newtheorem{rema}[theo]{Remark}

\def \kbar {{\bar k}}
\newcommand{\dem}[1] {\paragraph{Proof{#1}: }}

\def \Romannumeral #1 {\expandafter\uppercase\expandafter {\romannumeral #1} }
\def \br {{\rm{Br\,}}}
\def \sc{{\rm sc}}

\def \T {{\mathcal T}}
\def \C {{\mathcal C}}
\def \D {{\mathcal D}}
\def \XX {{\mathcal X}}
\def \HH {{\mathcal H}}

\def \nr {{\rm nr}}
\def \tors{{\rm tors}}

\def \pic {{\rm {Pic\,}}}
\def \upic {{\rm {UPic\,}}}
\def \div {{\rm{Div\,}}}
\def \gal {{\rm{Gal}}}

\def \calo {{\cal O}}
\def \spec {{\rm{Spec\,}}}

\def \Hom {{\rm {Hom}}}

\def \ext {{\rm {Ext}}}

\def \id {{\rm{id\,}}}

\def \Z {{\bf Z}}
\def \Q {{\bf Q}}

\def \F {{\bf F}}

\def \id {{\rm{id}}}
\def \G {{\bf G}_m}

\def\Sha{\cyrille X}

\def\ov{\overline}

\def\Ga{\Gamma}

\def\A{{\bf A}}

\def\smallsquare{\vbox{\hrule\hbox{\vrule height 1 ex\kern 1 ex\vrule}\hrule}}
\def\enddem{\hfill\smallsquare\vskip 3mm}

\def \abstract{\paragraph{Abstract. }}

\DeclareFontFamily{U}{wncy}{}
\DeclareFontShape{U}{wncy}{m}{n}{%
	<5>wncyr5%
		<6>wncyr6%
	<7>wncyr7%
	<8>wncyr8%
	<9>wncyr9%
	<10>wncyr10%
	<11>wncyr10%
	<12>wncyr6%
	<14>wncyr7%
	<17>wncyr8%
	<20>wncyr10%
	<25>wncyr10}{}
\DeclareMathAlphabet{\cyrille}{U}{wncy}{m}{n}
\def\Sha{\cyrille X}
\def\Be{\cyrille B}
\def\Ga{\Gamma}

\def \ss{{\rm ss}}
\def \sc{{\rm sc}}

\def \ab{{\rm ab}}

\def \tor{{\rm tor}}

\def\Br{{\rm Br\,}}

\def\Gm{{\mathbf{G}_m}}

\author{Cyril Demarche and David Harari}

\title{Local-global principles for homogeneous spaces of 
	reductive groups over global function fields}

\begin{document}
	\maketitle
	
\abstract{Let $K$ be a global field of positive characteristic. We prove that the Brauer-Manin obstructions to the Hasse principle, to weak approximation and to strong approximation are the only ones for homogeneous spaces of reductive groups with reductive stabilizers. The methods involve abelianization techniques and arithmetic duality theorems for complexes of tori over $K$.}
	
	\section{Introduction}\label{intro}
	
	Let $K$ be a global field of characteristic $p \geq 0$ and let $\A_K$
	denotes the ring of ad\`eles of $K$. Let $G$ be a reductive group 
	over $K$, and $X$ be a homogeneous space of $G$. We are interested in rational points on $X$, and more precisely, on various local-global principles associated to $X$: does $X$ satisfy the Hasse principle, i.e. does $X(\A_K) \neq \emptyset$ imply $X(K) \neq \emptyset$? If not, can we explain the failure using the so-called Brauer-Manin obstruction to the Hasse principle? Assuming that $X(K) \neq \emptyset$, can we estimate the size of $X(K)$ by studying the so-called weak and strong approximation on $X$ (with a Brauer-Manin obstruction if necessary), i.e. the closure of the set $X(K)$ in the topological space $X(\A_K^S)$, where $S$ is a (not necessarily finite) set of places of $K$ and $\A_K^S$ is the ring of $S$-ad\`eles (with no components in $S$)?
	
	The answer to those questions is known in the case where $K$ is a number field, provided that the stabilizers of points in $X$ are connected
	(see \cite{Borcrelle} and \cite{BDcom}).
	
	In the case of a global field of positive characteristic, the answer is known for semisimple simply connected groups (thanks to works by Harder, Kneser, Chernousov, Platonov, Prasad), but the general case is essentially open
	(see \cite{Ros-Tam}, Theorem 1.9 for some related results). In this paper, 
	we deal with these questions when both $G$ and the stabilizers are smooth, 
	connected and {\it reductive}. 
	
	\smallskip

Several new ingredients are needed to obtain our results: 

\begin{itemize}

\item To show that the Brauer-Manin obstruction to the Hasse principle 
is the only one (see Theorem~\ref{hptheo} for the precise result),
one has to use Poitou-Tate duality for complexes of tori in positive
characteristic (proven in \cite{DHdualcompl})
and a (non straightforward) compatibility result between Brauer-Manin and 
Poitou-Tate pairings. 

\item The statement on weak approximation (Theorem~\ref{weakaptheo})
relies on some part of 
Poitou-Tate exact sequence (which is established in \cite{DHdualcompl}) 
for a certain complex of tori, and on abelianization techniques 
(namely Lemma~\ref{compatweak}). Beforehand we define
in section~\ref{sectcompat} a new abelianization map associated 
to a homogeneous space $X=G/H$ as above, and prove a rather intricate
compatibility formula (Theorem~\ref{theo compatible strong approx}).
		
\item Theorem~\ref{mainstrong} presents the obstruction to strong approximation.
As in the number field case (settled in \cite{BDcom}), it is related to 
the Brauer-Manin pairing, but there are two important differences. The first 
one is that 
there is an additional term in the exact sequence describing the obstruction,
which reflects the fact that the global reciprocity map of class field theory 
is not surjective in positive characteristic. The second difference is that a
fibration method like in loc. cit. would probably not work here (see for instance section 5 in loc. cit.). 
Therefore, one should again rely on abelianization techniques
(in particular the compatibility formula of Theorem~\ref{theo compatible strong approx} in its full generality, and not only for elements of the 
algebraic Brauer group $\br_1 X$). An important role is also 
played by duality theorems for complex of tori, some of them extending 
results of \cite{DHdualcompl}.

\bigskip

		{\bf Notation and conventions.}

		\bigskip

\end{itemize}

	In the whole article (except in section~\ref{sectcompat},
	where $k$ is an arbitrary field), we consider 
	a finite field $k$ and a projective, smooth $k$-curve 
	$E$. We set $K = k(E)$,
	which is a global function field of characteristic $p$,
	and fix a separable closure $K^s$ of $K$. The absolute Galois group 
	$\gal(K^s/K)$ of $K$ is denoted by $\Ga_K$.
	Denote by $\Omega_K$ the set of all places of $K$; for every $v \in \Omega_K$,
	we will identify the Brauer group $\br K_v$ of the completion $K_v$ 
	to $\Q/\Z$ thanks to local class field theory. For every $K$-variety 
	$X$, we set $\ov X=X \times_K K^s$. The (cohomological) 
	Brauer group of $X$ is denoted $\br X$, and we set 
	$$\br_1 X:= \ker [\br X \to \br \ov X].$$
	We still denote by $\br K$ the image of $\br K$ in $\br X$, even though 
	the map $\br K \to \br X$ is not necessarily injective if $X$ has no 
	rational point. Notation like $H^i(K,C)$ for a commutative 
	$K$-group 
	scheme $C$ (resp. a bounded complex of commutative $K$-group schemes) always 
	denotes ffpf cohomology (resp. fppf hypercohohomology) of $C$. It coincides
	with \'etale (=Galois) cohomology when $C$ (resp. every group scheme occurring 
	in $C$) is smooth. For every finite set of places $S$ of $K$, we set 
	$$\Sha^i_S(K,C):=\ker [H^i(K,C) \to \prod_{v \not \in S} H^i(K_v,C)].$$
	$$\Sha^i(K,C)=\Sha^i_{\emptyset} (K,C); \quad \Sha^i_{\omega}(K,C)=
	\varinjlim_S \Sha^i_S(K,C),$$
	where the direct limit runs over all finite subsets $S$ of $\Omega_K$.
	The Pontryagin dual $A^D$ of a topological group $A$ is the group of 
	continuous homomorphism from $A$ to $\Q/\Z$ (if topology 
	is not specified, we assume that $A$ is discrete).
	
	\smallskip
	
	Let $G$ be a reductive group (always meaning: smooth connected reductive) 
	over $K$. Let $G^\ss$ denote the derived subgroup of $G$ and $G^\sc$
	the simply connected cover of $G^\ss$, together with the obvious morphism
	$\rho : G^\sc \to G$. Set $G^{\tor}=G/G^{\ss}$ (it is the maximal toric
	quotient of $G$). Let $T^\sc \subset G^\sc$ and $T_G \subset G$ be maximal
	tori such that $\rho(T^\sc) \subset T_G$. Let $C_G$ be the complex
	$C_G := [T^\sc \xrightarrow{\rho} T_G]$, with $T_G$ in degree $0$.
	Following Borovoi (cf. \cite{BorAMS} in characteristic zero),
	we have a natural map of Galois (hyper)cohomology sets:
	\[ \ab^1_G : H^1(K,G) \to H^1_{\ab}(K,G):=H^1(K,C_G) \, , \]
	which is functorial in $K$. There is an exact sequence of $K$-group 
	schemes 
	$$1 \to \mu_G \to G^{\sc} \to G^{\ss} \to 1,$$
	where $\mu_G$ is a finite $K$-group scheme of multiplicative type. It induces
	an exact triangle
	$$\mu_G[1] \to C_G \to G^{\tor} \to \mu_G[2].$$
	We denote by $Z_G$ the center of a reductive group $G$.
	
	\section{Hasse principle for homogeneous spaces} \label{two}
	
	We start with extending a well-known result on the abelianization maps to the 
	positive characteristic case. 
	
	\begin{prop} \label{prop ab1 bij} 
		There is a natural exact sequence of groups $$G^\sc(K) \to G(K) \xrightarrow{\ab^0_G} H^0(K, C_G) \to 1$$ and the map $\ab^1_G : H^1(K,G) \to H^1(K,C_G)$ is a bijection.
		
		For any place $v$ of $K$, there is an exact sequence $G^\sc(K_v) \to G(K_v) \xrightarrow{\ab^0_G} H^0(K, C_G) \to 1$ and the map
		$\ab^1_G : H^1(K_v,G) \to H^1(K_v,C_G)$ is a bijection.
	\end{prop}
	
	\dem{} Let $F$ be $K$ or $K_v$.
	
	By construction, one has a short exact sequence of groups and pointed sets
	\[G^\sc(F) \xrightarrow{\rho_*} G(F) \xrightarrow{\ab^0_G} H^0(F,C_G) \to H^1(F, G^\sc) \xrightarrow{\rho_*} H^1(F,G) \xrightarrow{\ab^1_G} H^1(F,C_G) \, .\]
	
	For $F=K$, following \cite{harder}, Satz A, the set $H^1(F, G^\sc)$ is trivial, hence the map $\ab^1_G$ has trivial kernel and the first sequence in the statement is exact. A twisting argument implies that the map $\ab^1_G$ is even injective. For the local statement ($F=K_v$), the injectivity (and exacteness of the first sequence) follows from \cite{BT}, Theorem 4.7.(ii).
	
	The proof of the surjectivity is an adaptation of the proofs of \cite{BorAMS}, Theorems 5.4 and 5.7, except that the existence of anisotropic maximal tori over local fields of positive characteristic is provided by \cite{DeB}, Lemma 2.4.1 (see also \cite{RosBig}, Proposition 4.4).
	\enddem

	We are interested in the Hasse principle for homogeneous spaces under $G$, with \emph{reductive} (recall that by definition this includes smoothness and 
	connectedness)
	geometric stabilizers. Following Raynaud (see \cite{Ray}, definition VI.1.1
	and Proposition VI.1.2), if $G$ is a smooth group scheme over $K$,
	a homogeneous space of $G$ is a smooth $K$-scheme $X$ with an action of $G$,
	such that for any $x \in X(K^s)$ (such an $x$ exists since $X$ is smooth),
	the stabilizer of $x$ in $\ov G:=G \times_K K^s$ is a finite type (over $K^s$) 
	subgroup scheme $H_x$ of $\ov G$
	and $\ov X$ is isomorphic to the quotient $\ov G/H_x$.
	
	\begin{defi}
		{\rm Let $X$ be a homogeneous space of a reductive group $G$ with 
			reductive stabilizer $\ov H=H_x$. 
			Let $L_X$ be the $K$-kernel defined by $X$ and $\sigma_X \in H^2(K, L_X)$ be the Springer class, that is the class 
			of the gerbe associated to $X$ (cf. \cite{FSS}, \S 5.2. or
			\cite{BorDuke}, \S 7.7  in characteristic zero). 
			By assumption, the stabilizers are reductive, hence following loc. cit.,
			there exists a $K$-torus $T$, which is a $K$-form of $\overline{H}^\tor$, and a natural map of marked sets
			\[ \ab^2_X : H^2(K,L_X) \to H^2(K,T) \, , \]
			which is functorial in $K$.
			On the other hand we have the class $\eta_X \in H^1(K,[T \to G^{\tor}])$
			constructed by Borovoi in \cite{borann} (where characteristic zero is
			assumed, but not used in the definition of $\eta_X$). There is an exact sequence 
			in Galois (hyper)-cohomology
			\begin{equation} \label{long}
			H^1(K,T) \to H^1(K,G^{\tor}) \to H^1(K,[T \to G^{\tor}]) \to 
			H^2(K,T) \to H^2(K,G^{\tor}).
			\end{equation}
		}
	\end{defi}
	
	\begin{prop} \label{twoclass}
		The image of $\eta_X$ in $H^2(K,T)$ is $\ab^2_X(\sigma_X)$.
	\end{prop} 
	
	\dem{} The method is similar to the one used in the proof of 
	\cite{BvH2}, Th. 9.6. We start with the case when $G$ itself is 
	a torus. Then $T$ is a subtorus of $G$ and $[T \to G]$ is quasi-isomorphic 
	to the quotient group $G/T$.
	Besides $X$ is a principal 
	homogeneous space of $G/T$ with class $[X] \in H^1(K,G/T) \simeq 
	H^1(K, [T \to G])$ corresponding to $\eta_X$. As $\sigma_X \in 
	H^2(K,T)$ is (by definition) just the image of $[X]$ by the 
	coboundary map $H^1(K,G/T) \to H^2(K,T)$, the result holds in this 
	case. 
	
	\smallskip 
	
	Assume now that $G$ is an arbitrary reductive group but satisfies 
	the additional hypothesis: 
	
	\smallskip
	
	(*) The canonical morphism $T \to G^{\tor}$ is injective. 
	
	\smallskip

	Then $Y:=X/G^{\ss}$
	is a homogenous space of $G^{\tor}$ with stabilizer (defined over $K$ in this case) 
	$T$. The Springer 
	class $\sigma_Y \in H^2(K,T)$ is (by construction) 
	just $\ab^2_X(\sigma_X)$, and it is 
	also the image of $\eta_Y \in H^1(K,[T \to G^{\tor}])$ in $H^2(K,T)$ 
	by the first case. But $\eta_Y=\eta_X$ by functoriality of the class 
	$\eta_X$ (\cite{borann}, \S 1.6), whence the result when (*) is 
	satisfied.
	
	\smallskip
	
	We now deal with the general case. By \cite{BvH2}, Prop~9.9 
	(whose proof is identical in characteristic $p$ thanks to the assumption
	that $G$ and $\ov H$ are reductive), there exists 
	a homogeneous space $Z$ of a $K$-group $F=G \times P$, where $P$ is 
	a quasi-trivial torus, such that: the homogeneous space $Z$ satisfies 
	(*) and there is a $K$-morphism $\pi : Z \to X$, compatible with the 
	respective actions of $F$, $G$ (via the projection $F \to G$), which makes 
	$Z$ an $X$-torsor under $P$. 
	In particular the $G$-homogeneous space $Z$ still has geometric stabilizer $\ov H$,
	with $L_Z=L_X$ and $\sigma_Z=\sigma_X$.
	Since Proposition~\ref{twoclass} holds for $Z$ 
	(because its satisfies (*)), it also holds for $X$ by 
	functoriality of the class $\eta_X$ and commutativity of the diagrams
	
	$$
	\begin{CD}
	H^2(K,L_Z) @>{\ab^2_Z}>> H^2(K, T \times P) \cr
	@VVV @VVV \cr
	H^2(K,L_X) @>{\ab^2_X}>> H^2(K,T) 
	\end{CD}
	$$
	
	\smallskip
	
	$$
	\begin{CD}
	H^1(K,[T,G^{\tor} \times P]) @>>> H^2(K, T \times P) \cr
	@VVV @VVV \cr
	H^1(K,[T,G^{\tor}]) @>>> H^2(K,T)
	\end{CD}
	$$
	
	\enddem
	
	\begin{rema}
		{\rm The previous proposition holds over an arbitrary field. 
			Recall also that the
			existence of a $K$-point on $X$ implies that the class $\sigma_X$ 
			is neutral, as well as the vanishing of $\eta_X$.}
		
	\end{rema}
	
	\begin{theo} \label{hptheo}
		Let $G$ be a reductive group over $K$. Then the Brauer-Manin obstruction to the Hasse principle is the only one for homogeneous spaces of $G$ with reductive stabilizers. More precisely, such a torsor $X$ has a rational point if and only if $X$ has an adelic point orthogonal to the subgroup $$\Be(X) := \ker[ (\br X/\br K) \to 
		\prod_{v \in \Omega_K} (\br X_{K_v}/\br K_v) \ ]$$
		for the Brauer-Manin pairing~: 
		\begin{equation} \label{bmp}
		X (\A_K) \times (\br X/\br K) 
		\to \Q/\Z, \quad (P_v) \mapsto
		\langle \alpha, (P_v) \rangle_{BM}:=
		\sum_{v \in \Omega_K} \alpha(P_v).
		\end{equation}
	\end{theo}
	
	Recall that by global class field theory, the sum 
	$\sum_{v \in \Omega_K} \alpha(P_v)$ is zero for every element $\alpha \in 
	{\rm Im} [\br K \to \br X]$, hence the Brauer-Manin pairing (\ref{bmp}) 
	is well defined. Also, for $\alpha \in \Be(X)$, the element 
	$\langle \alpha, (P_v) \rangle_{BM}$ is independent of the choice of 
	$(P_v) \in X (\A_K)$, because the localisation $\alpha_v \in 
	\br X_{K_v}$ of $\alpha$ is a constant element (i.e. it comes from 
	$\br K_v$) for every place $v$.

	\dem{}
	Fix a point $x \in X(K^s)$ and let $\ov H = H_x$
	be the stabilizer of $x$ in $\ov G$.
	
	Up to replacing $G$ by a flasque resolution (\cite{CTfl}) 
	$$1 \to S \to G' \to G \to 1,$$
	where $S$ is a flasque torus (which is central in $G'$) 
	and $G'$ is a quasi-trivial group 
	(that is: extension of a quasi-trivial torus by a semisimple
	simply connected group), we can assume that the group $G$ itself 
	is quasitrivial. Indeed $X$ is also a homogenous space of $G'$ such that 
	the stabilizer of $x$ is connected (it is an extension of $H_x$ by 
	$S_{K^s}$). In particular $\pic \ov G=0$ and the group of characters 
	$\widehat {G^{\tor}}$ of $G^{\tor}$ is a permutation Galois module
	(\cite{CTfl}, Prop~2.2).
	
	\smallskip
	
	Since $H^1(K,G^{\tor})=0$
	(the torus $G^{\tor}$ being quasi-trivial), we have by exact sequence 
	(\ref{long}) that the canonical map 
	$H^1(K,[T \to G^{\tor}]) \to H^2(K,T)$ is injective. Besides the class
	$\eta_X$ (viewed as an element of $H^2(K,T)$) is just $\ab^2_X(\sigma_X)$
	by Proposition~\ref{twoclass}.
	
	\smallskip
	
	Assume that $X(\A_K) \neq \emptyset$. Then the class $\sigma_X$ is neutral at every place $v$ of $K$, which implies that
	$\eta_X=\ab^2_X(\sigma_X) \in \Sha^2(K,T)$.
	The Brauer-Manin pairing defines a morphism $\beta_X : \Be(X) \to \Q / \Z$ 
	(recall that if $\alpha \in \Be(X)$, then $\langle \alpha, (P_v) \rangle_{BM} 
	\in \Q/\Z$ is independent of $(P_v) \in X(\A_K)$).
	
	\smallskip
	
	By \cite{BDH}, \S4, there is a complex $\upic(X)$
	(up to a shift this is 
	the complex ${\rm UPic} \ov X$ of \cite{BvH2}; we will recover this
	complex in section~\ref{sectcompat})
	such that $\br_1 X/\br K$ 
	is isomorphic to $H^1(K,\upic(X))$ (this is valid over any field $K$ 
	such that $H^3(K,\G)=0$).
	Moreover by \cite{BvH2}, Th. 5.8 (whose proof is 
	still valid in characteristic $p$ thanks to the additional assumption 
	that $G$ and $H$ are reductive), the complex $\upic(X)$
	is quasi-isomorphic 
	to $[\widehat{G^{\tor}} \to \widehat T]$ (recall that $\pic \ov G=0$).
	This induces a natural morphism $\lambda : \widehat T \to \upic(X)$ of $\Ga_K$-modules.
	Since $H^1(K,\widehat{G^{\tor}})=0$ (indeed $\widehat{G^{\tor}}$ is a permutation 
	Galois module), we obtain an injective morphism of abelian groups 
	$\lambda_* : H^1(K,\widehat T) \to \br_1 X/\br K$, which induces an 
	injective morphism of abelian groups
	$$\psi_X:=\lambda_* : \Sha^1(K,\widehat{T}) \hookrightarrow \Be(X).$$
	
	By Theorem 5.2 in \cite{DHdualcompl}, there is a Poitou-Tate perfect
	duality of finite groups
	\begin{equation} \label{duality sha2}
	\langle , \rangle_{PT} : \, 
	\Sha^2(K,T) \times \Sha^1(K,\widehat{T}) \to \Q / \Z \, . 
	\end{equation}
	This pairing and the class $\eta_X=\ab^2(\sigma_X)$
	define a morphism $\alpha_X :  \Sha^1(K,\widehat{T}) \to \Q/ \Z$.
	
	\begin{lem}
		The following diagram
		\[
		\xymatrix{
			\Sha^1(\widehat{T}) \ar[d]_{\psi_X} \ar[rd]^{\alpha_X} & \\
			\Be(X) \ar[r]^{\beta_X} & \Q / \Z
		}
		\]
		is commutative (up to sign).
	\end{lem}
	
	\dem{} By \cite{BvH2}, Th.~9.6. (whose proof is identical in characteristic 
	$p$ thanks to our additional assumptions),
	the class $\eta_X \in H^1(K,[T \to G^{\tor}]) \hookrightarrow 
	H^2(K,T)$ coincides (up to a sign) with the element of 
	$\ext^2(\upic(X),\G)$ given by the map $w : \upic(X) \to \G[2]$ of 
	exact sequence (1) in \cite{HS-op-desc}. Therefore the class 
	$\eta_X$ viewed in $\Sha^2(T) \subset H^2(K,T)=\ext^2(\widehat T,\G)$ is just 
	$\partial(\lambda)$, where $\partial : 
	\Hom_{\Ga_K}(\widehat T,\upic(X)) \to H^2(K,T)$ is the map defined in exact sequence (2)
	of \cite{HS-op-desc}. Now Theorem~3.5. of loc.cit. (whose proof
	in characteristic $p$ is identical) shows that for every $a \in
	\Sha^1(\widehat T)$, we have 
	$$\langle \eta_X ,a \rangle_{PT}=\langle \partial(\lambda),a \rangle_{PT}=
	\beta_X(\lambda_*(a)),$$
	or in other words~:
	$$\alpha_X(a)=\beta_X(\psi_X(a)).$$ This concludes the proof of the lemma.

	\enddem
	
	We can now finish the proof of Theorem~\ref{hptheo}.
	Assuming that $X(\A_K)^{\Be(X)} \neq \emptyset$, we have $\beta_X = 0$,
	hence $\alpha_X = 0$ by the lemma. The exactness of the
	pairing \eqref{duality sha2} and Proposition~\ref{twoclass} 
	imply that $\ab^2_X(\sigma_X)=\eta_X=0$.
	Now, the analogue of \cite{BorDuke}, Propostion 6.5 for global fields
	of positive characteristic (the proof of which is the same, except that
	Lemma 5.7 in \cite{BorDuke} is replaced by Proposition \ref{prop ab1 bij})
	implies that the map $\ab^2_X$ has "trivial kernel", hence $\eta_X$
	is neutral. As a consequence, there exists a $G$-equivariant
	morphism $P \to X$ defined over $K$, where $P$ is a $K$-torsor under $G$.
	Since $G$ is quasi-trivial, $H^1(K,G)$ is trivial (see \cite{harder}, Satz A),
	hence $P(K) \neq \emptyset$, thereore $X(K) \neq \emptyset$. 
	\enddem

\section{Abelianization of homogeneous spaces and a compatibility formula} \label{sectcompat}

Let $H$ be a reductive subgroup of a reductive group $G$ over an arbitrary 
field $k$. Set $X=G/H$ and denote by $e \in X(k)$ the image of the neutral element of $G$ in $X$.
Consider the complex of Galois-modules (or of commutative smooth $k$-group
schemes) defined
by $$C_X := [T_{H^\sc} \to T_H \oplus T_{G^\sc} \to T_G],$$ with $T_G$ in degree $0$.
In other words, we have $C_X := \textup{Cone}(C_H \to C_G)$. 
Denote by $\ab^0_X$ (see \cite{DemEdinburgh}, \S 2.6, where the characteristic zero
assumption is not needed thanks to 
our additional assumptions that $G$ and $H$ are smooth and 
reductive) the abelianization map $$\ab^0_X : X(k) \to H^0_{\ab}(k,X):=H^0(k,C_X).$$ 
Following \cite{DemJOA}, we set 
$$\br_1(X,G)=\ker [\br X \to \br \ov G],$$
which is a subgroup of $\br X$ containing $\br_1 X$. We also consider 
the subgroup $\Br_{1,e}(X,G) \simeq \br_1(X,G)/\br k$ of $\br_1(X,G)$
consisting 
of those elements $\alpha$ such that $\alpha(e)=0$, and use a similar notation 
for $\br_{1,e} X \simeq \br_1 X/ \br k$. 

\smallskip

The goal of this section is to prove
(see Theorem~\ref{theo compatible strong approx} below) 
that there is a natural isomorphism 
\[\phi_X :  H^1(k, \widehat{C}_X) \xrightarrow{\sim} \Br_{1,e}(X,G) \, , \]
and that the following diagram
\[
\xymatrix{
	X(k) \ar@<-30pt>[d]_{\ab'_X} \times \Br_{1,e}(X,G) \ar[r]^(.65){\textup{ev}} & \Br(k) \ar[d]^= \\
	H^0(k, C_X) \times H^1(k, \widehat{C}_X) \ar@<-30pt>[u]_{\phi_X} \ar[r]^(.7){\cup} & \Br(k)
}
\]
is commutative ($\widehat{C}_X$ denotes the dual complex of $C_X$: it
is the cone of the morphism of complexes
\[ [\widehat{T}_G \to \widehat{T}_{G^\sc}] \to [\widehat{T}_H \to \widehat{T}_{H^\sc}] \, \]
where $\widehat{T}_H$ is in degree $0$). Here 
$\ab'_X$ is an abelianization map, which we are going to define 
by changing a little the map $\ab^0_X$. 
In \cite{Demeio}, Demeio proves such a compatibility in characteristic zero 
for $\ab^0_X$ with a long and impressive cocyle computation. It is 
quite likely that the maps $\ab^0_X$ and our new map $\ab'_X$ actually 
coincide, but we did not succeed in proving this. The modified 
map $\ab'_X$ seems more suitable to check the required compatibility. 

\subsection{Abelianization map over a field} \label{subsec ab field}

From now on, let $k$ be a field, $X$ a smooth geometrically integral $k$-variety and $\pi : Y \to X$ a torsor under a reductive $k$-group $H$. We fix a point $y_0 \in Y(k)$ and let $x_0 := \pi(y_0) \in X(k)$. Following \cite{DemJOA}, we define $Z$ to be $Y/Z_H$, $z_0$ to be the image of $y_0$ and $\upic'(\pi)$ to be the following complex of Galois modules:
\[\upic'(\pi) := [\kbar(Z)^\times \to \div(\overline{Z}) \to \pic'(\overline{Z}/\overline{X})] \, ,\]
with $\kbar(Z)^\times$ in degree $-1$ (here $\pic'(\overline{Z}/\overline{X})$
is the relative Picard group of $\ov Z$ over $\overline{X}$).
We also define
\[\upic(\pi) := [\kbar(Z)^\times/\kbar^\times \to \div(\overline{Z}) \to \pic'(\overline{Z}/\overline{X})] \, ,\]
with the obvious natural exact sequence of complexes:
\[0 \to \kbar^\times[1] \to \upic'(\pi) \to \upic(\pi) \to 0 \, .\]

We will sometimes need the following pointed version of those complexes, which are canonically quasi-isomorphic to the previous ones:
\[\upic'(\pi)_0 := [\kbar(Z)_0^\times \to \div(\overline{Z})_0 \to \pic'(\overline{Z}/\overline{X})]\]
and
\[\upic(\pi)_0 := [\kbar(Z)_{0,1}^\times \to \div(\overline{Z})_0 \to \pic'(\overline{Z}/\overline{X})] \, ,\]
where $\kbar(Z)_0^\times$ (resp. $\kbar(Z)_{0,1}^\times$) denotes the subgroup of rational functions defined at $z_0$ (resp. taking the value $1$ at $z_0$) and $\div(\overline{Z})_0$ is the group of divisors $D$ such that $z_0$ is not contained in the support of $D$. We also have the classical complexes 
$\upic(Z)=[\kbar(Z)^*/\kbar^* \to \div \ov Z]$, $\upic'(Z)$ etc.
(corresponding to the case when $\pi : Z\to Z$ is the identity map).

The construction of $\upic'(\pi)$ is contravariant in $\pi$ (see for instance Proposition 2.2 in \cite{DemJOA}), and for any $x \in X(k)$, the natural morphism $\kbar^\times[1] \to \upic'(\pi_x : Z_x \to \spec k)$ is a quasi-isomorphism.

Therefore, one gets a well-defined map
\[{\ab''}_\pi : X(k) \to \Hom_k(\upic'(\pi), \kbar^\times[1]) \, ,\]
where $\Hom_k$ denotes the set of morphisms in the derived category of bounded complexes of Galois modules.

In addition, when a point $y_0 \in Y(k)$ is given, one gets a natural splitting $\upic(\pi) \to \upic'(\pi)$, hence it defines the required map
\[{\ab'}_\pi : X(k) \to \Hom_k(\upic(\pi), \kbar^\times[1]) \, .\]
When $Y=X$ and $\pi = \id_X$, we denote ${\ab'}_\pi$ by ${\ab'}_X$.

We want to think about ${\ab'}_\pi$ as an abelianization map for the set of rational points of $X$, which is a replacement of $\ab^0_X$ when $X=G/H$.

We now study the d\'evissage of $\upic(\pi)$ in terms of $Y$ and $H$ :
\begin{lem} \label{lem triangle exact}
	There is a natural exact triangle 
\[\upic(\pi) \to \upic(Y) \to \upic(H) \to \upic(\pi)[1] \, .\]
\end{lem}

\dem{}
This is essentially the proof of Corollary 3.3 in \cite{DemJOA}. Let 
$H'=H/Z_H$ be the quotient of $H$ by its center. 
Since $\upic(\pi)_0[1]$ is the cone of $\upic(Z)_0 \to \upic(H')_0$, 
the following commutative diagram of complexes
\[\xymatrix{
\upic(Z)_0 \ar[r] \ar[d] & \upic(H')_0 \ar[d] \\
\upic(Y)_0 \ar[r]^{\varphi} & \upic(H)_0 
}
\]
induces a canonical morphism of complexes
\[\alpha : \upic(\pi)_0[1] \to \textup{cone}(\varphi) \, ,\]
such that the following diagram of exact triangles is commutative :
\[
\xymatrix{
\upic(Z)_0 \ar[r] \ar[d] & \upic(H')_0 \ar[d] \ar[r] & \upic(\pi)_0[1] \ar[d]^\alpha \ar[r] & \upic(Z)_0[1] \ar[d] \\
\upic(Y)_0 \ar[r]^{\varphi} & \upic(H)_0 \ar[r] & \textup{cone}(\varphi) \ar[r] & \upic(Y)_0[1] \, .
}\]

Let us now prove that $\alpha$ is a quasi-isomorphism. Since the complexes $\upic(\pi)_0[1]$ and $\textup{cone}(\varphi)$ are concentrated in degrees $-2$ to $0$, we only compute the cohomology corresponding to those degrees : one has a commutative diagram of long exact sequences
\begin{changemargin}{-3cm}{2cm}
\begin{equation} \label{diag qis}
\xymatrix{
0 \ar[r] & H^{-2}(\upic(\pi)[1]) \ar[r] \ar[d]^{\alpha^{-2}} & U(\overline{Z}) \ar[r] \ar[d] & \widehat{H'} \ar[r] \ar[d] & H^{-1}(\upic(\pi)[1]) \ar[r] \ar[d]^{\alpha^{-1}} & \pic(\overline{Z}) \ar[r] \ar[d] & \pic(\overline{H'}) \ar[r] \ar[d] & H^{0}(\upic(\pi)[1]) \ar[r] \ar[d]^{\alpha^0} & 0 \\
0 \ar[r] & H^{-2}(\textup{cone}(\varphi)) \ar[r] & U(\overline{Y}) \ar[r] & \widehat{H} \ar[r]^(.3){\beta} & H^{-1}(\textup{cone}(\varphi)) \ar[r] & \pic(\overline{Y}) \ar[r] & \pic(\overline{H}) \ar[r] & H^{0}(\textup{cone}(\varphi)) \ar[r] & 0 \, .
}
\end{equation}
\end{changemargin}
The proof of Corollary 3.3 in \cite{DemJOA} ensures that $\alpha^{-2}$ and $\alpha^0$ are isomorphisms.

Let us now prove that $\alpha^{-1}$ is an isomorphism (the proof of this fact in \cite{DemJOA} is too sketchy): using the exact sequences $\widehat{Z_H} \to \pic(\overline{Z}) \to \pic(\overline{Y})$ and $\widehat{Z_H} \to \pic(\overline{H'}) \to \pic(\overline{H})$, diagram chasing in \eqref{diag qis} proves that the result follows from the commutativity of the following diagram:
\[
\xymatrix{
\pic(\overline{X}) & H^{-1}(\upic(\pi)_0[1]) \ar[l]_(.6){\sim} \ar[d]^{\alpha^{-1}} \\
\widehat{H} \ar[r]^(.3){\beta} \ar[u]^\Delta & H^{-1}(\textup{cone}(\varphi)) \, ,
}
\]
where $\Delta : \widehat{H} \to \pic(\overline{X})$ is the map defined by $\chi \mapsto \chi_* [\overline{Y}]$ and the isomorphism $$H^{-1}(\upic(\pi)_0[1]) \to \pic(\overline{X})$$ is constructed in \cite{DemJOA}. To prove the required commutativity, given $\chi \in \widehat{H}$, functoriality of the various maps implies that it is sufficient to consider the case $H = \Gm$ (and $H'=1$, $Z=X$) and $\chi = \id : \Gm \to \Gm$. More precisely, let $D \in \div(\overline{X})$ with support not containing $x_0$, let $\pi : Y_D \to \overline{X}$ be the associated $\Gm$-torsor. Let $U \subset \overline{X}$ be the complement of the support of $D$. By construction, the pullback $Y_{D,U} \to U$ of $Y_D \to \overline{X}$ admits a canonical section $s$ over $U$, inducing a map $f : Y_{D,U} \to \Gm$ which can be seen as an element $f \in \overline{k}(Y_D)^\times$. Let $y_0 := s(x_0) \in Y_{D,U}(k)$ and let $Y_0$ be the fiber of $\pi$ at $x_0$. Then the restriction of $f$ at $Y_0$ and the point $y_0 \in Y_0(k)$ induce a morphism $f_0 : \Gm \to \Gm$. Making explicit the maps $\alpha^{-1}$ and $\beta$, the required commutativity boils down to the natural equalities $\pi^*(D) = \textup{div}(f)$ in $\div(\overline{Y_D})$ and $f_0 = \id$.
\enddem

We now want to compare the map ${\ab'}_\pi$ defined earlier with the maps $\ab^0_H$ and $\ab^1_H$ defined by Borovoi in \cite{BorAMS} for reductive $k$-groups $H$. We first prove the following

\begin{lem} \label{lem compat ab ab'}
With the above notation, we have a commutative diagram (up to a sign) with exact rows
\begin{changemargin}{-3.5cm}{2cm}
\[
\xymatrix{
H(k) \ar@{^{(}->}[r] \ar[d]^{\ab^0_H} & Y(k) \ar[d]^{{\ab'}_Y} \ar[r]^\pi &  X(k) \ar[r] \ar[d]^{{\ab'}_\pi} & H^1(k,H) \ar[d]^{\ab^1_H} \\
 H^0_{\ab}(k,H) = \Hom_k(\upic(H), \overline{k}^\times[1]) \ar[r] & \Hom_k(\upic(Y), \overline{k}^\times[1]) \ar[r]^{\pi'} & \Hom_k(\upic(\pi), \overline{k}^\times[1]) \ar[r] & H^1_{\ab}(k,H) = \Hom_k(\upic(H), \overline{k}^\times[2]) \, , 
}\]
\end{changemargin}
where the second line comes from Lemma \ref{lem triangle exact} and the isomorphisms $$H^i_{\ab}(k,H) = \Hom_k(\upic(H), \overline{k}^\times[i+1])$$ come from \cite{BvH}, Theorem 4.8 and Corollary 4.9.

In addition, when $Y=G$ is a reductive $k$-group, $H \subset G$ is a reductive subgroup and $X := G/H$, then the map ${\ab'}_G : G(k) \to \Hom_k(\upic(G), \kbar^\times[1])=H^0_\ab(k,G)$ coincides with $\ab^0_G$.
\end{lem}

\dem{}
The commutativity of the central square is a consequence of the functoriality of the map ${\ab'}_\pi$ with respect to morphisms of torsors.

Let us now prove the commutativity of the left hand side square : using the functoriality of the map ${\ab_\pi}$ and the fact that the morphism $H \to Y$, given by the action on $y_0 \in Y(k)$, is $H$-equivariant, it is sufficient to prove that ${\ab^0}_H = {\ab'}_H$, where ${\ab'}_H : H(K) \to \Hom_k(\upic(H), \overline{k}^\times[1])$ is defined as ${\ab'}_{\id_H}$. 
There exists a coflasque resolution (see \cite{CTfl}, Proposition 4.1)
\[1 \to P \to H_1 \to H \to 1 \, ,\]
where $P$ is a quasi-trivial $k$-torus and $H_1$ is an extension of
a (coflasque) 
torus $T$ by a semisimple simply connected $k$-group $H^\ss$. Then functoriality of the map ${\ab'}_H$ and Hilbert 90 imply that it is enough to prove the required compatibility ${\ab'}_T = {\ab^0}_T$ for the $k$-torus $T$. By definition, the map ${\ab'}_T : T(k) \to \Hom_k(\widehat{T}[1], \kbar^\times[1])=\Hom_k(\widehat{T}, \kbar^\times)$ is given by $t \mapsto (\chi \mapsto \chi(t))$, and it clearly coincides with the map ${\ab^0}_T: T(k) \to T(k)$, composed with the natural identification $T(k) \xrightarrow{\sim} \Hom_k(\widehat{T}, \kbar^\times)$. It concludes the proof of the commutativity of the left hand side square. Note that this also proves the last statement in Lemma \ref{lem compat ab ab'}.

Let us now prove the remaining commutativity, concerning the right hand side square. Let $x \in X(k)$ and consider the torsor $\pi_x : Y_x \to \spec k$ defined as the pullback of $\pi$ by $x$ ($Y_x$ denotes the fiber of $Y$ at $x$). We have a commutative diagram 
\[
\xymatrix{
\upic(\pi) \ar[d]^{\alpha} \ar[r]^{y} & \upic'(\pi) \ar[d]^{\alpha'} \ar[r]^{x^*} & \upic'(\pi_x) \ar[d]^{\alpha_x} & \kbar^\times[1] \ar[l] \ar[d]^= \\
\textup{cone}(\varphi)[-1] \ar[r]^y & \textup{cone}(\varphi')[-1] \ar[r]^{x^*} & \textup{cone}(\varphi'_x)[-1] & \kbar^\times[1] \ar[l] \, , 
}
\]
where the vertical maps are quasi-isomorphisms (see the proof of Lemma \ref{lem triangle exact}).
We now prove that the diagram
\[
\xymatrix{
\textup{cone}(\varphi)[-1] \ar[r]^y & \textup{cone}(\varphi')[-1] \ar[r]^{x^*} & \textup{cone}(\varphi'_x)[-1] & \kbar^\times[1] \ar[l] \\
\upic(H)[-1] \ar[u] \ar[r]^= & \upic(H)[-1] \ar[u] \ar[r]^= &  \upic(H)[-1] \ar[u] & \upic(Y_x)[-1] \ar[u]^{\partial_{Y_x}} \ar[l]^{\varphi_x}
}
\]
is commutative. The only non obvious square is the right hand side one, but it comes from the functoriality of the cone in the category of complexes, together with the definition of $\partial_{Y_x}$.

The two previous commutative diagrams imply that the square
\[
\xymatrix{
Y(k) \ar[r] \ar[d]^{{\ab'}_\pi} & H^1(k,H) \ar[d]^{{\ab'}^1_H} \\
\Hom_k(\upic(\pi), \overline{k}^\times[1]) \ar[r] & \Hom_k(\upic(H), \overline{k}^\times[2])
}
\]
commutes, where ${\ab'}^1_H$ maps the cohomology class of 
a torsor $W \to \spec k$ under $H$ to the rightmost morphism in the natural exact triangle
\[\kbar^\times[1] \to \upic'(W) \to \upic(W) \xrightarrow{\partial_W} \kbar^\times[2]\] 
with the canonical identification $\upic(W) \xrightarrow{\sim} \upic(H)$.

We conclude the proof using \cite{BvH}, Theorem 5.5, which proves that ${\ab'}^1_H = - \ab^1_H$.
\enddem

We need to prove other properties of the map $\ab'_\pi$ :
\begin{prop} \label{prop ab action}
	Let $X=G/H$ be the quotient of a reductive group $G$ by a reductive subgroup $H$. Let $\pi : G \to X$ be the quotient map (pointed by $e \in G(k)$).
	Then for all $g \in G(k)$ and $x \in X(k)$, 
	\[{\ab'}_\pi(g \cdot x) = \pi'(\ab^0_G(g)) + {\ab'}_\pi(x) \, . \]
	In particular, if $G$ is semi-simple and simply connected, we have 
	${\ab'}_\pi(g \cdot x) = {\ab'}_\pi(x) \, $.
\end{prop}

\dem{}
Consider $\id_G : G \to G$ as a torsor under the trivial group, and $\id \times \pi : G \times_k G \to G \times_k X$ as a natural torsor under $H$. Then we have natural morphisms of torsors :
\[
\xymatrix{
G \ar[d]^{\id_G} & G \times_k G \ar[r]^{p_2} \ar[l]^{p_1} \ar[d]^{\id_G \times \pi} & G \ar[d]^\pi \\
G & G \times_k X \ar[l]^{p_G} \ar[r]^{p_X} & X \, ,
}
\]
that induce, by functoriality of $\upic(\pi)$ and by Lemma \ref{lem triangle exact}, a commutative diagram in the derived category, where the rows are exact triangles (see \cite{BvH}, Lemma 5.1 for the third vertical map) :
\[
\xymatrix{
\upic(H)[-1] \ar[r] \ar[d]^\sim & \upic(G) \oplus \upic(\pi) \ar[r] \ar[d]^{p_G^* + p_X^*} & \upic(G) \oplus \upic(G) \ar[d]^\sim \ar[r] & \upic(H) \ar[d]^\sim \\
\upic(1 \times H)[-1] \ar[r] & \upic(\id_G \times \pi) \ar[r] & \upic(G \times_k G) \ar[r] & \upic(1 \times H) \, .
}
\]
The five-lemma implies that the norphism $p_G^* + p_X^* : \upic(G) \oplus \upic(\pi) \to \upic(\id_G \times \pi)$ is an isomorphism.

By functoriality of the construction of ${\ab'}_\pi$, the morphism of torsors 
\[
\xymatrix{
G \times_k G \ar[r]^{m_G} \ar[d]^{\id_G \times \pi} & G \ar[d]^\pi \\
G \times_k X \ar[r]^m & X
}
\]
induces a commutative diagram
\[
\xymatrix{
& \kbar^\times[1] & \kbar^\times[1] \oplus \kbar^\times[1] \ar[l]_{\times} \\
& \upic(\id_G \times \pi) \ar[u]_{(g,x)} & \upic(G) \oplus \upic(\pi) \ar@<0.8cm>[u]^g \ar@<-0.8cm>[u]^{x} \ar[l]_{p_G^* + p_X^*} \\
\upic(\pi) \ar[ruu]^{g \cdot x} \ar[ru]^{m^*} \ar[rru]_{(\textup{can}, \id)} & & \, ,} 
\]
which concludes the proof. The non trivial commutativity in this last diagram is that of the triangle at the bottom, which we explain now: recall that we are given a point $x_0 = \pi(e) \in X(k)$. Consider the following commutative diagram of morphisms of torsors:
\[
\xymatrix{
G \ar[r]^(.3){(\id_G,1)} \ar[d]^{\id_G} & G \times G \ar[d]^{(\id_G, \pi)} & G \ar[d]^{\pi}  \ar[l]_(.3){(1, \id_G)} \\
G \ar[r]^{\iota_{x_0}} & G \times X & X \ar[l]_{\iota_1} \, ,
}
\]
where the bottom horizontal maps are defined by $\iota_{x_0}(g) = (g, x_0)$ and $\iota_1(x) = (1,x)$. If we denote by $\varpi_G$ (resp. $\varpi_X$) the projection from $\upic(G) \oplus \upic(\pi)$ to $\upic(G)$ (resp. $\upic(\pi)$), then we deduce from the previous diagram that $\varpi_G = \iota_{x_0}^* \circ (p_G^* + p_X^*)$ and $\varpi_X = \iota_1^* \circ (p_G^* + p_X^*)$. But $m \circ \iota_{x_0} = \pi : G \to X$ and $m \circ \iota_1 = \id_X : X \to X$, therefore we get that the required triangle commutes.

\smallskip

Finally, if $G$ is assumed to be semi-simple and simply connected, 
then the complex $\upic (G)$ is quasi-isomorphic to $0$, hence 
the map $\pi'$ is trivial, which implies the required result.

\enddem

\begin{rema}
	{\rm
	It is worth noting that the construction of the map $\ab_\pi$ depends on the choice of a $k$-point $y_0 \in Y(k)$. But one can prove, using the same kind of arguments as in the proof of Lemma \ref{prop ab action}, that the map $\ab_\pi$ depends only on the image $x_0$ of $y_0$ in $X(k)$. More precisely, two points $y_0, y_0' \in Y(k)$ such that $\pi(y_0) = \pi(y_0')$ define the same map $\ab_\pi$, or equivalently, the construction of this map depends only on the choice of a point $x_0 \in \pi(Y(k))$.
	}
\end{rema}

Lemma \ref{lem compat ab ab'} and Proposition~\ref{prop ab action} 
imply the following proposition :

\begin{prop} \label{orbitesprop}
	If $k$ is a non archimedian local field and $X = G/H$ is the quotient of a reductive group $G$ by a reductive subgroup $H$, then the map $\ab'_\pi : X(k) \to \Hom_k(\upic(\pi), \kbar^\times[1])$, associated to the torsor $\pi : G \to X$, is surjective, and for all $x, x' \in X(k)$, $\ab'_\pi(x) = \ab'_\pi(x')$ if and only if there exists $g \in G^\sc(k)$ such that $x' = g \cdot x$. 
	The same holds when $k$ is a global field of positive characteristic.
\end{prop}

\dem{} In the proof, $k$ denotes either a non archimedian local field, or a global field of positive characteristic.

Consider the following commutative diagram with exact rows (see Lemma \ref{lem compat ab ab'}):
\[
\xymatrix{
H(k) \ar[r] \ar[d]^{\ab^0_H} & G(k) \ar[r]^\pi \ar[d]^{\ab^0_G} & X(k) \ar[r] \ar[d]^{\ab'_\pi} & H^1(k,H) \ar[d]^{\ab^1_H} \ar[r] & H^1(k,G) \ar[d]^{\ab^1_G} \\
H^0_\ab(k,H) \ar[r] & H^0_\ab(k,G) \ar[r]^(.35){\pi'} & \Hom_k(\upic(\pi), \kbar^\times [1]) \ar[r] & H^1_\ab(k,H) \ar[r] & H^1_\ab(k,G) \, .
}
\]
Then diagram chasing proves that surjectivity of $\ab'_\pi$ is a consequence of surjectivity of $\ab^1_H$, injectivity of $\ab^1_G$, surjectivity of $\ab^0_G$ (those properties follows from Proposition \ref{prop ab1 bij}), and of Proposition \ref{prop ab action}.

Let us now prove the second part of the Proposition: let $x$, $x' \in X(k)$ such that $\ab'_\pi(x) = \ab'_\pi(x')$. Using the previous diagram and the injectivity of $\ab^1_H$ (see Proposition \ref{prop ab1 bij}), we get that $x$ and $x'$ have the same image in $H^1(k,H)$. Therefore, there exists $g \in G(k)$ such that $g \cdot x = x'$. Applying Lemma \ref{prop ab action}, we get that $\pi'(\ab^0_G(g)) = 0$, hence by Proposition \ref{prop ab1 bij}, $g$ lifts to $G^\sc(k)$, which concludes the proof.
\enddem

Let $X=G/H$, with $G$ and $H$ reductive. Let us now construct a canonical isomorphism $\phi_X : \widehat{C}_X \to \upic(\pi)$ in the derived category, inspired by \cite{DemJOA}, sections 4.1.2 and 4.1.3.

By construction, $\widehat{C}_X$ is the cone of the morphism of complexes 
\[ [\widehat{T}_G \to \widehat{T}_{G^\sc}] \to [\widehat{T}_H \to \widehat{T}_{H^\sc}] \, \]
where $\widehat{T}_H$ is in degree $0$.
Consider the following commutative diagram of complexes, where the vertical maps are either obvious or defined in \cite{BvH}, section 4):
\begin{equation} \label{diag upic char}
\xymatrix{
 [\widehat{T}_G \to \widehat{T}_{G^\sc}] \ar[r] \ar[d] &  [\widehat{T}_H \to \widehat{T}_{H^\sc}] \ar[d] \\
\textup{Cone}(\upic(T_G)_0 \to \upic(T_{G^\sc})_0) \ar[r] & \textup{Cone}(\upic(T_H)_0 \to \upic(T_{H^\sc})_0) \\
\textup{Cone}(\upic(G)_0 \to \upic(G^\sc)_0) \ar[u] \ar[r] \ar[d] \ar[u] & \textup{Cone}(\upic(H)_0 \to \upic(H^\sc)_0) \ar[u] \ar[d] \\
\upic(G)_0 \ar[r]^\varphi & \upic(H)_0 \\
\upic(Z)_0 \ar[u] \ar[r] & \upic(H')_0 \ar[u] \, . 
}
\end{equation}
All the vertical maps, except the ones between the last two lines, are quasi-isomorphisms. Hence this diagram induces a natural isomorphism $\phi'_X : \widehat{C}_X \to \textup{Cone}(\varphi)$, which we can compose with the quasi-isomorphism $\alpha : \upic(\pi)[1] \to \textup{Cone}(\varphi)$ induced by the two last lines of the previous diagram (see the proof of Lemma \ref{lem triangle exact}), to get a natural isomorphism $\phi_X : \widehat{C}_X \to \upic(\pi)[1]$.

By construction, this isomorphism fits into the following commutative diagram of exact triangles in the derived category:
\begin{equation} \label{diag devissage C_X}
\xymatrix{
\widehat{C}_H[-1] \ar[r] \ar[d]^{\phi_H} & \widehat{C}_X \ar[r] \ar[d]^{\phi_X} & \widehat{C}_G \ar[r] \ar[d]^{\phi_G} & \widehat{C}_H \ar[d]^{\phi_H} \\
\upic(H)[-1] \ar[r] & \upic(\pi) \ar[r] & \upic(G) \ar[r] & \upic(H) \, ,
}
\end{equation}
where the vertical maps are isomorphisms (morphisms $\phi_H$ and $\phi_G$ are defined in \cite{BvH}, section 4, and also in the first four lines of the previous big diagram).

In addition, Lemma 3.1 in \cite{BDH} implies that the natural morphism $H^0(k, C_X) \to \Hom_k(\widehat{C}_X, \kbar^\times[1])$ is an isomorphism. These facts lead to the following:
\begin{defi} \label{def ab'}
We denote by ${\ab'}_X : X(k) \to H^0(k, C_X)$ the map defined by 
\[X(k) \xrightarrow{{\ab'}_\pi} \Hom_k(\upic(\pi), \kbar^\times[1]) \xrightarrow{\widehat{\phi_X}} \Hom_k(\widehat{C}_X, \kbar^\times[1]) \xleftarrow{\sim} H^0(k, C_X) \, ,\]
where $\pi : G \to X$ is the quotient morphism.
\end{defi}

In particular, Lemma \ref{lem compat ab ab'}, Definition \ref{def ab'} and diagram \eqref{diag devissage C_X} imply that the following useful diagram is commutative (up to sign):
\begin{equation} \label{diag devissage ab}
\xymatrix{
	H(k) \ar[r] \ar[d]^{\ab^0_H} & G(k) \ar[r]^\pi \ar[d]^{\ab^0_G} & X(k) \ar[r] \ar[d]^{{\ab'}_X} & H^1(k,H) \ar[r] \ar[d]^{\ab^1_H} & H^1(k,G) \ar[d]^{\ab^1_G} \\
	H^0(k,C_H) \ar[r] & H^0(k,C_G) \ar[r] & H^0(k,C_X) \ar[r] & H^1(k,C_H) \ar[r] & H^1(k,C_G) \, ,
}
\end{equation}
where the unnamed maps are the natural ones.

\subsection{The compatibility result}

We can now prove the main compatibility result of this section, which can be seen as a generalization of \cite{DemJOA}, Theorem 4.14 and of \cite{BDH}, Theorem 6.2:

\begin{theo} \label{theo compatible strong approx}
Let $k$ be a field and $G$ be a reductive group.
Let $H \subset G$ be a reductive $k$-subgroup and $X := G/H$. Let $\pi : 
	G \to X$ be the quotient map (pointed by $e \in G(k)$ and its
	image $x_0:=\pi(G)$).
The canonical isomorphism $\phi_X : \widehat{C}_X \to \upic(\pi)$ in the derived category induces an isomorphism 
\[\phi_X :  H^1(k, \widehat{C}_X) \xrightarrow{\sim} \Br_{1,e}(X,G) \, , \]
and the following diagram
\[
\xymatrix{
X(k) \ar@<-30pt>[d]_{{\ab'}_X} \times \Br_{1,e}(X,G) \ar[r]^(.65){\textup{ev}} & \Br(k) \ar[d]^= \\
H^0(k, C_X) \times H^1(k, \widehat{C}_X) \ar@<-30pt>[u]_{\phi_X} \ar[r]^(.7){\cup} & \Br(k)
}
\]
is commutative, up to a universal sign (independent of all the data).
\end{theo}

\dem{}
Following \cite{DemJOA}, there is a natural morphism $\upic'(\pi) \to \tau_{\leq 2} R \pi_* \Gm[1]$ inducing an isomorphism $H^1(k, \upic(\pi)) \xrightarrow{\sim} \Br_{1,e}(X,G)$. Together with the isomorphism $\phi_X : \widehat{C}_X \to \upic(\pi)$, we get the required isomorphism. Let us now prove the commutativity.

For any $x \in X(k)$, one can define a natural splitting $\ab'(x)$ of $\kbar^\times \to \tau_{\leq 2} R \pi_* \Gm$ induced by $x$.

We can decompose the diagram above as the composition of the following diagrams:
\[
\xymatrix{
X(k) \ar@<-30pt>@/_10pc/[ddd]_{{\ab'}_\pi} \ar@<-30pt>[d]_{{\ab'}} \times \Br_{1,e}(X,G) \ar[r]^(.6){\textup{ev}} & \Br(k) \ar[d]^= \\
\Hom_k(\tau_{\leq 2} R \pi_* \Gm [1], \kbar^\times[1])  \ar@<-30pt>[d] \times H^1(k, \tau_{\leq 2} R \pi_* \Gm [1]) \ar@<-30pt>[u] \ar[r]^(.8)\cup &  \Br(k) \ar[d]^= \\
\Hom_k(\upic'(\pi), \kbar^\times[1]) \ar@<-30pt>[d]_e \times H^1(k, \upic'(\pi)) \ar@<-30pt>[u] \ar[r]^(.75)\cup &  \Br(k) \ar[d]^= \\
\Hom_k(\upic(\pi), \kbar^\times[1]) \ar@<-30pt>[d]_{\widehat{\phi_X}} \times H^1(k, \upic(\pi)) \ar@<-30pt>[u]_e \ar[r]^(.75)\cup &  \Br(k) \ar[d]^= \\
\Hom_k(\widehat{C}_X, \kbar^\times[1])  \ar@<-30pt>[d]_\sim \times H^1(k, \widehat{C}_X) \ar@<-30pt>[u]_{\phi_X} \ar[r]^(.7){\cup} & \Br(k) \ar[d]^= \\
H^0(k, C_X) \times H^1(k, \widehat{C}_X) \ar@<-30pt>[u]_= \ar[r]^(.6){\cup} & \Br(k) \, .
}
\]
By construction, and by functoriality of cup-products, this last diagram is commutative (up to sign).
\enddem

\begin{rema}
	{\rm
	Let $X = G/H$ as in Definition \ref{def ab'} and Theorem \ref{theo compatible strong approx}. Up to replacing $G$ by a flasque resolution $G_1$,
	one can realise $X$ as the quotient of the quasi-trivial group $G_1$
	by a reductive subgroup $H_1$. Since $\pic(\overline{G_1})=0$, one gets a natural isomorphism 
	$$\widehat{C_X} \xrightarrow{\sim} [\widehat{G_1} \to \widehat{T_{H_1}} \to \widehat{T_{H_1^\sc}}],$$ where $\widehat{G_1}$ is a permutation
	Galois module. Assuming the group $G$ is quasi-trivial
	will be very useful in the next two sections.}
\end{rema}

\subsection{The abelianization map over an arbitrary base}
In this section, we extend the definition of the map $\ab'_X$ for homogeneous spaces of reductive group schemes defined over an arbitrary base scheme $S$. It will be useful in the next sections in order to take integral points into account (case $S = \spec(\mathcal{O}_v)$).

Let $S$ be an integral regular noetherian scheme and $H$ be a reductive group scheme over $S$. Let $\pi : Y \to X$ be a torsor under $H$. For any $S$-scheme $W$, let $p_W : W \to S$ denote the structure morphism. 

Let $Z := Y / Z_H$ and $\varpi : Z \to X$ be the associated $H' := H/Z_H$-torsor. We define $\upic'(\pi)$ to be the following complex of \'etale sheaves over $S$ :
\[\underline{\upic'}(\pi) := [{p_Z}_* \underline{K}_{Z/X}^\times \to {p_Z}_* \underline{\div}_{Z/X} \to {p_X}_* \underline{\pic}_{Z/X}] \, ,\]
where the sheaves $\underline{K}_{Z/X}^\times$ and $\underline{\div}_{Z/X}$ are defined in \cite{HS-op-desc}, Appendix A. By loc. cit., there is a natural exact sequence of \'etale sheaves over $X$
\[0 \to (\Gm)_X \to \varpi_* \underline{K}_{Z/X}^\times \to \varpi_* \underline{\div}_{Z/X} \to \underline{\pic}_{Z/X} \to 0 \, .\] 
Applying ${p_X}_*$, one gets a natural morphism $(\Gm)_S \to \underline{\upic'}(\pi)$, and we define $\underline{\upic}(\pi)$ to be the cone of this morphism, whence an exact triangle
\[(\Gm)_S \to \underline{\upic'}(\pi) \to \underline{\upic}(\pi) \to (\Gm)_S[1] \, . \]

Let $y_0 \in Y(S)$ be a point, and let $x_0 := \pi(y_0)$. Denote 
by $D(S)$ the derived category of bounded complexes of étale sheaves over
$S$. Following the construction in section \ref{subsec ab field}, we get a natural map, functorial in $S$ and $X$,
\[\ab'_\pi : X(S) \to \Hom_{D(S)}(\underline{\upic}(\pi), (\Gm)_S) \, .\]

Let now $G$ be a reductive group scheme over $S$ and $H \subset G$ a reductive subgroup scheme. Taking $Y=G$ and $X=G/H$, we can apply the previous constructions. From now on, we assume that $H$ and $G$ admit compatible maximal tori $T_H \subset T_G$. The diagram \eqref{diag upic char}, as a diagram of \'etale sheaves over $S$, still holds. All vertical morphisms, except the bottom ones, are quasi-isomorphisms of \'etale sheaves over $S$, since it is true over any separably closed field (see after diagram \eqref{diag upic char}). Similarly, since the result holds over separably closed fields, diagram \eqref{diag upic char} induces an isomorphism of complexes of \'etale sheaves $\phi_X : \widehat{C}_X \to \underline{\upic}(\pi)[1]$.

As a conclusion, composing the map $\ab'_\pi : X(S) \to \Hom_{D(S)}(\underline{\upic}(\pi), (\Gm)_S)$ with the isomorphism $\phi_X$, one gets the required map:
\[\ab'_X : X(S) \to H^0(S, C_X) \, ,\]
that is functorial in $S$ and $(H,G)$, and that coincides with the definition of section \ref{subsec ab field} in the case $S$ is the spectrum of a field.

\section{Weak approximation} \label{four}

From now on, the setting is the following: $K=k(E)$ is again 
the function field of a projective, smooth curve $E$ over a finite 
field $k$. We consider a reductive linear algebraic group $G$
over $K$, and $X$ is a homogeneous space of $G$. We assume that
$X(K) \neq \emptyset$. Let $e \in X(K)$ and let $H \subset G$ be the stabilizer of $e$ in $G$. Then $X = G/H$ (with $e$ identified to the image of 
the neutral element of $G$ in $X$) and we still suppose that 
$H$ is a reductive subgroup of $G$. Set $X(K_{\Omega})=\prod_{v \in \Omega_K} 
X(K_v)$. 

\smallskip

We are interested in the closure
of $X(K)$ in $X(K_\Omega)$, for the product topology. We define 
$\Be_\omega(X)$  
as the subgroup of $\br_{1,e} \, X \simeq \br_1 X/\br K$
consisting of those
elements $\alpha$ such that their localization
$\alpha_v \in \br_{1,e} \, X_{K_v}$ is zero for all but finitely many $v$.
For each finite set of places $S$ of $K$, we set $K_S=\prod_{v \in S} K_v$
and we define $$\Be_S(X):=\ker [\br_{1,e} \, X  \to \prod_{v \not \in S} 
\br_{1,e} \, X_{K_v}].$$  
In particular $\Be(X)$ (cf. section~\ref{two}) 
identifies to $\Be_{\emptyset}(X)$. 
For all $\alpha \in \Be_\omega(X)$, $(P_v) \in X(K_{\Omega})$, 
the Brauer-Manin pairing~: 
$$\langle \alpha,(P_v) \rangle_{BM}=\sum_{v \in \Omega_K} \alpha(P_v)$$
is well-defined.  
For every subgroup $B$ of $\Be_\omega(X)$, we denote
by $X(K_{\Omega})^B \subset X(K_{\Omega})$ the orthogonal of $B$ for
the Brauer-Manin pairing.

\begin{theo} \label{weakaptheo}
	Let $G$ be a reductive group over $K$ and $H$ a reductive subgroup.
	Then the Brauer-Manin obstruction to weak approximation on $X=G/H$
	associated to $\Be_\omega(X)$ is the only one, i.e. $X(K)$ is dense
	in $X(K_\Omega)^{\Be_\omega(X)}$.
	More precisely, for any finite set $S$ of places of $K$,
	the Brauer-Manin pairing induces a surjective map
	$X(K_S) \to \left(\Be_S(X) / \Be(X) \right)^D$, whose kernel is exactly
	the closure of $X(K)$ inside $X(K_S)$.
\end{theo}

\begin{rema}
	{\rm Assume that $X$ admits a regular compactification, i.e. there exists a regular proper $K$-variety $X^c$ and an open immersion $X \to X^c$. Then the group $\Be_\omega(X)$ is exactly the algebraic Brauer group of $X^c$ (see for
example \cite{BDH}, Prop~4.1).}
\end{rema}

\dem{} 
Up to replacing $G$ by a flasque resolution, 
one can assume that $G$ is quasi-trivial.
Using an elementary instance of the fibration method, we see
that $G$ satisfies weak approximation, since $G^\sc$ (by \cite{prasad}) 
and the quasi-trivial (hence $K$-rational) 
torus $G^\tor$ do satisfy weak approximation.
In addition, we have $H^1(K,G) = 1$ by \cite{harder} and Hilbert's 90 
(the latter shows that $H^1(K,G^{\tor})=0$ thanks to Shapiro's lemma).

\smallskip

Let $C := [H^\tor \to G^\tor]$ and
$C_X := \textup{Cone}(C_H \to C_G)$ (cf. section~\ref{sectcompat}).
We have the 
Cartier duals $\widehat C=[\widehat {G^\tor} \to \widehat {H^\tor}]$ and 
$\widehat {C_X}=\textup{Cone}(\widehat {C_G} \to \widehat {C_H})[-1]$.

\smallskip

By construction, we have a natural commutative diagram of exact triangles of complexes:
\begin{equation} \label{diag hom spaces devissage}
\xymatrix{
	\mu_H[1] \ar[d] & & \mu_H[2] \ar[r]^\sim \ar[d] & \mu_H[2] \ar[d] \\
	C_H \ar[r] \ar[d] & C_G \ar[r] \ar[d]^{\textup{qis}} & C_X \ar[r] \ar[d] & C_H[1] \ar[d] \\
	H^\tor \ar[r] \ar[d] & G^\tor \ar[r] & C \ar[r] \ar[d] & H^\tor[1] \ar[d] \\
	\mu_H[2] & & \mu_H[3] \ar[r]^\sim & \mu_H[3] \, .
}
\end{equation}

Since $G$ is assumed to be quasi-trivial, we have $G^{\ss}=G^{\sc}$, hence
$\mu_G=0$ and $C_G$ is quasi-isomorphic to $G^{\tor}$, which is a
quasi-trivial torus. Taking Cartier duals, we get an exact triangle~:
$$ \widehat {C_X} \to \widehat {G^{\tor}}[1]
\to \widehat {C_H} \to \widehat {C_X} [1],$$
whence (using $H^1(K,\widehat {G^{\tor}})=0$) an exact sequence~:
$$0 \to H^0(K,\widehat {C_H}) \to H^1(K,\widehat {C_X})
\to H^2(K,\widehat {G^{\tor}})$$
and similarly replacing $K$ with a completion $K_v$.
As $\widehat {G^{\tor}}$ is a permutation Galois module, we have 
$\Sha^2_S(K,\widehat {G^{\tor}})=0$ by Shapiro's lemma and \v Cebotarev's 
Theorem. Thus
$$\Sha^0_S(K,\widehat{C_H}) \xrightarrow{\sim} \Sha^1_S(K, \widehat{C_X}).$$
We also have an exact triangle 
$$  \widehat {\mu_H} [-2] \to 
\widehat C \to \widehat {C_X} \to \widehat {\mu_H} [-1],$$
which yields an isomorphism $$\Sha^1_S(K, \widehat{C}) \to 
\Sha^1_S(K, \widehat{C_X}).$$
Summing up, we get isomorphisms 
\[\Sha^0_S(K,\widehat{C_H}) \xrightarrow{\sim} \Sha^1_S(K, \widehat{C_X}) \xleftarrow{\sim} \Sha^1_S(K, \widehat{C}) \, .\] 

\smallskip

As recalled before (cf. proof of Th. 2.5.), the group $\br_{1,e} X \simeq 
\br_1 X/ \br K$
is isomorphic to $H^1(K,\upic(X))=H^1(K,\widehat C)$ (and this is true 
over any field). Therefore $\Sha^1_S(K,\widehat C)$ (resp.
$\Sha^1(K,\widehat C)$) identifies to $\Be_S(X)$ (resp. to 
$\Be(X)$), whence a Brauer-Manin pairing~: 
$$X(K_S) \times \Sha^1_S(K,\widehat C) \to \Q/\Z, ((P_v)_{v \in S}, \alpha)
\mapsto \sum_{v \in S} \alpha(P_v),$$
which is trivial on $X(K_S) \times \Sha^1(K,\widehat C)$. It induces a map 
${\rm BM} : X(K_S) \to 
(\Sha^1_S(K,\widehat{C}) / \Sha^1(K, \widehat{C}))^D$. On the other hand,
local duality for complex of tori (\cite{Demtores}, Th. 3.1) induces 
a map
$$\theta_S : \prod_{v \in S} H^1_{\ab}(K_v,H)=\prod_{v \in S} 
H^1(K_v, C_H) \to
\left(\Sha^0_S(K, \widehat{C_H})/\Sha^0(K, \widehat{C_H})\right)^D.$$ 
We also have a map in the derived category of Galois modules
$$C_X \otimes_{\bf L} \widehat{C_X} \to \G[1],$$ 
which induces (for ecah completion $K_v$) a cup-product pairing 
$$H^0(K_v,C_X) \times H^1(K_v,\widehat{C_X}) \to H^2(K_v,\G)=\Q/\Z.$$
The latter is compatible (in an obvious sense) with local duality 
for the complexes $C_H$, $\widehat{C_H}$ via the maps $C_X \to C_H[1]$
and $\widehat{C_H} \to \widehat{C_X}[1]$. Finally, denote by 
$\partial_K : X(K) \to H^1(K,H)$ the coboundary map,
which (by composing with the abelianization
map) yields a morphism of pointed sets
$\partial_K ^{\ab} : X(K) \to H^1_{\ab}(K,H)=H^1(K,C_H)$. Similarly
we have (for every completion $K_v$) maps
$\partial_v^{\ab} : X(K_v) \to H^1(K_v,C_H)$.
By Lemma~\ref{lem compat ab ab'}, the 
map $\partial_v^{\ab}$ can also be obtained by composing the 
abelianization map 
$\ab ' _{X_{K_v}} \to H^0(K_v,C_X)$ (defined in section~\ref{sectcompat}) 
with the natural map $H^0(K_v,C_X) \to H^1(K_v,C_H)$ (induced by 
diagram (\ref{diag hom spaces devissage})). 

\begin{lem} \label{compatweak}
	There is a commutative diagram with exact rows~:
	\begin{equation} \label{diag esp hom WA}
	\xymatrix{
		G(K) \ar[r] \ar[d] & X(K) \ar[r]^{\partial_K^{\ab}}
		\ar[d] & H^1(K,C_H) \ar[r] \ar[d] & 1 \\
		G(K_S) \ar[r] & X(K_S) \ar[d]^{BM} \ar[r]^{\partial_S^{\ab}=\prod_{v \in S} \partial_v^{\ab} } &
		\prod_{v \in S} 
		H^1(K_v,C_H) \ar[r] \ar[d]^{\theta_S} & 1 \\
		& \left(\Sha^1_S(K,\widehat{C}) /
		\Sha^1(K, \widehat{C})\right)^D \ar[r]^{\sim}
		& \left(\Sha^0_S(K, \widehat{C_H})/\Sha^0(K, \widehat{C_H})\right)^D \ar[d] & \\
		& & 0 & \, .
	}
	\end{equation}
\end{lem}

\dem{} The exactness of rows follows from the triviality of $H^1(K,G)$ 
and $H^1(K_v,G)$ (recall that $G$ is quasi-trivial) combined with 
Proposition~\ref{prop ab1 bij}. The only non-trivial remaining point 
is the commutativity of the bottom square. By functoriality of 
the cup-product and Theorem~\ref{theo compatible strong approx},
there is a commutative diagram~:
$$
\begin{CD}
X(K_v) @>{(\ab ')_{X_{K_v}} }>> H^0(K_v,C_X) @>>> H^1(K_v,C_H) \cr
@VVV @VVV @VVV \cr
\br_{1,e}(X_{K_v},G_{K_v})^D @>{(\Phi_{X_{K_v}})^D}>>
H^1(K_v,\widehat{C_X})^D @>>> H^0(K_v,\widehat C_H)^D
\end{CD}
$$

where the left vertical map is given by the local evaluation pairing 
$$X(K_v) \times \br_{1,e}(X_{K_v},G_{K_v}) \to \Q/\Z, (P_v,\alpha_v) 
\mapsto \alpha_v(P_v)$$
and the right vertical map by local duality for the complex $C_H$.
Now let $a \in \Sha^0_S(K,\widehat{C_H})$ 
with image $b \in \Sha^1_S(K,\widehat{C_X})$ and $c \in \Sha^1_S(K,\widehat{C})$.
Set $\alpha=\Phi_X(b) \in \br_{1,e} X \subset \br_{1,e}(X,G)$.
Let $v \in S$ and $P_v \in X(K_v)$. 
By the previous diagram, we have 
$$\alpha(P_v) =(\partial_v^{\ab}(P_v) \cup a_v) \in \br K_v=\Q/\Z ,$$
where $a_v \in H^0(K_v,\widehat{C_H})$ is the localization of $a$.
Hence 
$$(BM((P_v))).c:=\sum_{v \in S} \alpha(P_v)=
\sum_{v \in S} (\partial_v^{\ab}(P_v) \cup a_v)=(\theta_S(\partial_S^{\ab}(P_v))).a,$$
which yields the required commutativity.

\enddem

\begin{rema}
	{\rm Actually we used the difficult compatibility proven in
	Theorem~\ref{theo compatible strong approx} only for those 
	elements of $H^1(K_v,\widehat C_X)$ 
	coming from $H^1(K_v,\widehat C)$, so Theorem~6.2. of 
	\cite{BDH} (whose proof is much easier) would 
	be sufficient at this stage. However, we will definitely need
	Theorem~\ref{theo compatible strong approx} in its full generality 
	in section~\ref{strongsect}.
	}
\end{rema}

We resume the proof of Theorem~\ref{weakaptheo}.
Let us now prove that the right hand side column of diagram
(\ref{diag esp hom WA}) is exact. Consider the following commutative diagram:
\[
\xymatrix{
	& \bigoplus_{v \notin S} H^1(K_v,C_H) \ar[r]^\sim \ar[d] & \left(\prod_{v \notin S} H^0(K_v, \widehat{C_H}) \right)^D \ar[d] & \\
	H^1(K, C_H) \ar[r] \ar[d]^= & \bigoplus_{v \in \Omega} H^1(K_v,C_H) \ar[r] \ar[d]^p & \left(H^0(K, \widehat{C_H}) / \Sha^0(K, \widehat{C_H})\right)^D \ar[r] \ar[d] & 0 \\
	H^1(K, C_H) \ar[r] & \prod_{v \in S} H^1(K_v,C_H) \ar[r] \ar[d] & \left(\Sha^0_S(K, \widehat{C_H}) / \Sha^0(K, \widehat{C_H})\right)^D \ar[r] \ar[d] & 0 \\
	& 0 & 0 & \, .
}
\]
Using \cite{DHdualcompl}, Theorem 5.7, the second row is exact. And by construction, the columns are exact. Hence an easy diagram chase implies that the bottom row is exact. Therefore, the right hand side column in \eqref{diag esp hom WA} is exact. In addition, we know that $G(K)$ is dense in $G(K_S)$. Therefore, an easy diagram chase in \eqref{diag esp hom WA}, together with the comparison Theorem 6.2 in \cite{BDH},  implies that the map 
\[
X(K_S) \xrightarrow{BM} \left(\Sha^1_S(K,\widehat{C}) / \Sha^1(K, \widehat{C})\right)^D = \left(\Be_S(X) / \Be(X) \right)^D \]
is surjective, and that the inverse image of $0$ is exactly the closure of $X(K)$, which concludes the proof.
\enddem

Theorem~\ref{weakaptheo} can be slightly refined when the homogeneous
space $X$ is a reductive group:

\begin{cor} \label{groupcase}

	Let $L$ be a reductive group over $K$. Let $C_L=[T^{\sc} \to T_L]$
	be the complex of tori associated to $L$.
	Then there is an exact sequence of groups

$$
1 \to \overline{L(K)} \to L(K_\Omega) \to \Sha^1_\omega(K, \widehat{C_L})^D \to \Sha^1(K,C_L) \to 1 \, .
$$

\end{cor}

\dem{} We can view $L$ as a homogeneous space $L=G/H$, with $G$ semi-simple
and simply connected, and $H$ reductive. By Theorem~\ref{weakaptheo}, there
is an exact sequence of pointed sets 
$$1 \to \overline{L(K)} \to L(K_\Omega) \stackrel{BM}{\to} \Be_{\omega}(L)^D \simeq \Sha^1_{\omega}(K, \widehat{C_L})^D.$$
By Lemma~\ref{compatweak}, the Brauer-Manin map $L(K_\Omega) \to
\Sha^1_{\omega}(K, \widehat{C_L})^D$ is the composition of the 
abelianization map $\ab_L : 
L(K_\Omega) \to \prod_{v \in \Omega_K} H^0(K_v,C_L)$ 
with the map (which is induced by local duality) $\theta : \prod_{v \in \Omega_K} H^0(K_v,C_L) \to \Sha^1_{\omega}(K, \widehat{C_L})^D$. Therefore
the Brauer-Manin map is a morphism of groups. By \cite{BorAMS}
(paragraph~3.10), 
there is (for every completion $K_v$) an exact sequence of pointed sets
$$L(K_v) \stackrel{\ab_L ^v} \to H^0(K_v,C_L) \to H^1(K_v,L^{\sc}),$$
which implies that $\ab_L$ is surjective because $H^1(K_v,L^{\sc})$ is trivial
(\cite{harder}, Satz A).

\smallskip

It is now sufficient to show that the sequence 
of abelian groups
$$\prod_{v \in \Omega_K} H^0(K_v,C_L) \to \Sha^1_\omega(K, \widehat{C_L})^D \to \Sha^1(K,C_L) \to 0$$
is exact. This is done by observing that by \cite{Demtores}, Th.~3.1. and 
\cite{DHdualcompl}, Th. 5.2, this sequence is the dual
of the exact sequence of discrete abelian groups 
$$0 \to \Sha^1(K,\widehat{C_L}) \to \Sha^1_\omega(K, \widehat{C_L}) \to
\bigoplus_{v \in \Omega_K} H^1(K_v, \widehat{C_L}).$$

\enddem

\begin{rema}
{\rm The results and proofs of sections \ref{two}, \ref{sectcompat} and 
\ref{four} are still valid 
over a number field. The analogues of 
Theorems~\ref{hptheo} and \ref{weakaptheo} were previously known in this
context: they are proven by Borovoi in \cite{Borcrelle} via more geometric 
techniques (namely fibration methods); the case of principal homogeneous
spaces is due to Sansuc \cite{San}. Another approach (over an 
arbitrary global field) is to use flasque resolutions, see 
\cite{nguyen}, Th.~3.9. 
In the next section, we will see that
the situation is slightly different for strong approximation.
}
\end{rema}

\section{Strong approximation} \label{strongsect}

Notation is as in the previous section.
We set $A^\star :=\Hom(A,\Z)$ for every abelian group $A$. 
An abelian group is said to be {\it in the class $\mathcal{E}$}
(cf. Definition 3.10 in \cite{DHdualcompl}) if it is an extension
of a finitely generated group by a profinite group.

\smallskip

Let $C=[T_1 \to T_2]$ be a short complex of $K$-tori with dual 
$\widehat C=[\widehat {T_2} \to \widehat{T_1}]$. For 
every $\chi \in H^{-1}(K,\widehat{C})$, denote by $\chi_v \in 
H^{-1}(K_v,\widehat{C})$ the localization of $\chi$ at the place $v$.
The cup-product pairing
$$H^0(K_v,C) \times H^{-1}(K_v,\widehat{C}) \to K_v^*,$$
induces a pairing 

\begin{equation} \label{pairing char c}
\begin{array}{cl}
H^0(\A_K, C) \times H^{-1}(K,\widehat{C}) & \to \Z \\
((g_v), \chi) & \mapsto \sum_v v(\chi_v.g_v) \cdot [k(v):k] \, .
\end{array}
\end{equation}
Since the degree of a principal divisor on the projective curve
$E$ is zero, this pairing is trivial on the subgroup 
$H^0(K,C) \times H^{-1}(K,\widehat{C})$.

\begin{lem} \label{wedgelem} 
	
	a) The kernel $P$ of the map $H^0(\A_K, C)/H^0(K,C) \to
	H^{-1}(K,\widehat{C})^\star$ induced by (\ref{pairing char c}) is 
	profinite. 
	
	\smallskip
	
	b) The canonical morphism $i : H^0(\A_K, C)_{\wedge}/H^0(K,C)_{\wedge} \to
	\left(H^0(\A_K, C)/H^0(K,C)\right)_\wedge$ is an isomorphism.
	
	\smallskip
	
	c) The map $H^0(\A_K,C) \to H^{-1}(K, \widehat{C})^\star$ is surjective and 
	the diagram
	\[
	\xymatrix{
		H^0(\A_K,C)/H^0(K,C) \ar[r] \ar[d] & \left(H^0(\A_K,C) / H^0(K,C) \right)_\wedge \ar[d] \\
		H^{-1}(K, \widehat{C})^\star \ar[r] & {H^{-1}(K, \widehat{C})^\star}_\wedge 
	}
	\]
	is cartesian.

\end{lem} 

\dem{}

a) Up to replacing $C$ by a quasi-isomorphic complex,
one can assume that $T_1$ is quasi-trivial. This yields a commutative
diagram with exact rows and surjective left vertical map:
$$
\begin{CD}
T_1(\A_K) @>>> T_2(\A_K) @>>> H^0(\A_K,C) @>>> 0 \cr
@VVV @VVV @VVV \cr
\widehat T_1(K)^{\star} @>>> \widehat T_2(K)^{\star} @>>> 
H^{-1}(K,\widehat C)^{\star} 
\end{CD}
$$

The snake lemma now implies that the kernel of
$H^0(\A_K,C)/H^0(K,C) \to H^{-1}(K, \widehat{C})^\star$ is a quotient
of $\ker\left(T_2(\A_K)/T_2(K) \to \widehat{T_2}(K)^\star\right)$
by a closed subgroup, and this kernel is profinite by
\cite{RosBig}, Proposition 5.7.5. Hence $\ker\left(H^0(\A_K,C)/H^0(K,C)
\to H^{-1}(K, \widehat{C})^\star \right)$ is profinite. 

\smallskip

b) The exact sequence 
$$H^0(K,C) \to H^0(\A_K,C) \to H^0(\A_K,C)/H^0(K,C) \to 0$$
and Lemma~3.12 a) of \cite{DHdualcompl}  show that $i$ is surjective. 
Since $H^{-1}(K, \widehat{C})^\star$ is a lattice, a) shows that 
$H^0(\A_K,C)/H^0(K,C)$ is in the class $\mathcal{E}$.
In particular, the
canonical map $H^0(\A_K,C)/H^0(K,C) \to (H^0(\A_K,C)/H^0(K,C))_{\wedge}$
is injective, which shows that $P$ injects into
$(H^0(\A_K,C)/H^0(K,C))_{\wedge}$ as well as in
$H^0(\A_K, C)_{\wedge}/H^0(K,C)_{\wedge}$. The commutative diagram 
\begin{equation} \label{diagwedge}
\begin{CD}
0 @>>> P @>>> H^0(\A_K,C)/H^0(K,C) @>>> H^{-1}(K,\widehat{C})^\star \cr
&& 	@V{=}VV @VVV @VVV \cr
0 @>>> P @>>> H^0(\A_K,C)_{\wedge}/H^0(K,C)_{\wedge}@>>>
H^{-1}(K,\widehat{C})^\star _{\wedge} \cr 
&&      @V{=}VV @VViV @VVV \cr
0 @>>> P @>>> (H^0(\A_K,C)/H^0(K,C))_{\wedge}@>>>
H^{-1}(K,\widehat{C})^\star _{\wedge} \cr
\end{CD}
\end{equation}
has exact first line by definition, and exact third line because 
it is obtained by completing the first line and
$H^{-1}(K,\widehat{C})^\star$ is a lattice,
so \cite{DHdualcompl}, Lemma~3.12 b) applies. To prove 
the injectivity of $i$, it is sufficient (by diagram chasing) to show
that the second line is exact as well. Let $\pi : H^0(\A_K,C) \to 
H^0(\A_K,C)/H^0(K,C)$ be the projection. The exact sequence 
$$0 \to \pi^{-1}(P) \to H^0(\A_K,C) \to H^{-1}(K,\widehat{C})^\star$$
induces (again by \cite{DHdualcompl}, Lemma~3.12 b) an exact sequence 
$$\pi^{-1}(P)_{\wedge} \to H^0(\A_K,C)_{\wedge}
\to H^{-1}(K,\widehat{C})^\star _{\wedge},$$
whence an exact sequence 
$$\pi^{-1}(P)_{\wedge}/H^0(K,C)_{\wedge} \stackrel{j}{\to}
H^0(\A_K,C)_{\wedge}/H^0(K,C)_{\wedge} \to
H^{-1}(K,\widehat{C})^\star _{\wedge}.$$
But the exact sequence 
$$H^0(K,C) \to \pi^{-1}(P) \to P \to 0$$
induces (as $P$ is profinite) 
a surjective map $u : \pi^{-1}(P)_{\wedge}/H^0(K,C)_{\wedge} \to P$
and $j$ factorizes through $u$, hence the second line of the diagram is 
exact, as required.

\smallskip 

c) If $C = [0 \to \Gm]$, the required surjectivity is obviously true.
If $C=[0 \to T]$, then one can find a resolution $0 \to R \to Q \to T \to 0$
(e.g. a flasque resolution) of $T$ by $K$-tori $R$ and $Q$,
with $Q$ quasi-trivial. It induces an injective map of lattices
$\widehat T(K) \to \widehat Q(K)$, hence a surjective map
$\widehat{Q}(K)^\star \to \widehat{T}(K)^\star$; thus
one can reduce to the case of a quasi-trivial torus.
By Shapiro's Lemma, this case reduces to the known case of
$\Gm$. Hence the map $H^0(\A_K,C) \to H^{-1}(K, \widehat{C})^\star$ is
surjective as soon as $C = [0 \to T]$ for any torus $T$.
Let $C = [T_1 \to T_2]$ be an arbitrary complex of tori.
Then the map $\widehat{T_2}(K)^\star \to H^{-1}(K, \widehat{C})^\star$
is surjective (again by injectivity of $H^{-1}(K,\widehat C)
\to \widehat C(K)$), hence the surjectivity of $H^0(\A_K,C) \to
H^{-1}(K, \widehat{C})^\star$ follows from that of
$T_2(\A_K) \to \widehat{T_2}(K)^\star$. Since the diagram (\ref{diagwedge}) 
is commutative with exact rows, the second point follows.

\enddem

\begin{prop} \label{complexc}
	a) There is an exact sequence
	\begin{equation} \label{exact seq complex}
	H^0(K,C) \to H^0(\A_K, C) \to \left(H^0(\A_K, C)/H^0(K,C)\right)_\wedge \xrightarrow{\partial} {H^{-1}(K,\widehat{C})^\star}_\wedge / H^{-1}(K,\widehat{C})^\star \to 0 \, ,
	\end{equation}
	where the morphism $\partial$ is given (after completion) by
	pairing~(\ref{pairing char c}) for $C$. 
	
	\smallskip
	
	b) There is an exact sequence 
	
	$$H^0(K,C) \to H^0(\A_K,C) \to (H^1(K,\widehat C)/\Sha^1(K,\widehat C))^D \stackrel{\partial_C}{\to} 
	{H^{-1}(K,\widehat{C})^\star}_\wedge / H^{-1}(K,\widehat{C})^\star \to 0.$$
\end{prop}

Observe that for $C=\G$, we have $\widehat C=\Z[1]$ and b) is just the 
classical exact sequence of global class field theory.

\dem{} By Theorems 5.7 and 5.10 in \cite{DHdualcompl}, there is a 
commutative diagram with exact rows~:
\begin{equation} \label{ptdiag}
\begin{CD}
H^0(K,C) @>>> H^0(\A_K,C) @>>>
(H^1(K,\widehat C)/\Sha^1(K,\widehat C))^D \cr
@VVV @VVV @VV=V \cr 
H^0(K,C)_{\wedge} @>>> H^0(\A_K,C)_{\wedge} @>>>
(H^1(K,\widehat C)/\Sha^1(K,\widehat C))^D @>>> 0  \quad , \cr 
\end{CD}
\end{equation}
whence exactness of the sequence $H^0(K,C) \to H^0(\A_K, C) \to
\left(H^0(\A_K, C)/H^0(K,C)\right)_\wedge$ thanks to 
Lemma~\ref{wedgelem} b). The exactness of sequence \eqref{exact seq complex} 
now follows from Lemma~\ref{wedgelem}, c), and part b) of the proposition
follows from its part a) and diagram (\ref{ptdiag}).

\enddem 

We now consider a homogenous space $X=G/H$ with $G$ and $H$ reductive and 
we assume further that $G$ is quasi-trivial. As in the previous sections, we
define $C=[H^{\tor} \to G^{\tor}]$ and $C_X={\rm Cone}(C_H \to C_G)=
[T_{H^{\sc}} \to T_H \to G^{\tor}]$. We denote by $\widehat C$ and 
$\widehat {C_X}$ their respective duals (cf. section~\ref{four}).
We have the analogue of the pairing (\ref{pairing char c}) with 
$C$ (resp. $\widehat C$) replaced by $C_X$ (resp. $\widehat{C_X}$), 
and again the pairing is trivial on
$H^0(K,C_X) \times H^{-1}(K,\widehat{C_X})$.

\smallskip

We observe that $H^{-1}(K,\widehat{C})=H^{-1}(K,\widehat{C_X})$ (and similarly 
over every completion $K_v$ of $K$) thanks to the exact triangle 
$$\widehat{\mu_H}[-2] \to \widehat C \to \widehat{C_X} \to \widehat{\mu_H}[-1].$$
For every finite and non empty set of places $S$ of $K$, we set
$$\br_S \, X=\ker [\br_e X \to \prod_{v \in S} \br X_{K_v} ]$$
(not to be confused with the groups $\Be_S(X)$ of section~\ref{four}) 
and $\br_S(X,G):=\br_S \, X \cap \br_{1,e}(X,G)$, 
$\br_{1,S} X:=\br_S \, X \cap \br_{1,e} X $.
We will also use a smooth model $\mathcal G$ (resp. $\mathcal H$, $\mathcal X=
\mathcal G/\mathcal H$)
of $G$ (resp. $H$, $X$) over some non empty Zariski open subset $U \neq E$ of 
the curve $E$. Shrinking $U$ if necessary, 
we can assume that $\HH$ and $\mathcal G$ admit compatible maximal tori 
$T_{\HH} \subset T_{\mathcal G}$.
We have the corresponding complexes $C_{\HH}$, 
$\mathcal C=[\mathcal H^{\tor} \to \mathcal G^{\tor}]$, 
$C_{\XX}$ defined over $U$. For every bounded complex $\mathcal F$ of flat 
commutative finite type group schemes over $U$, the 
compact support hypercohomology groups $H^i_c(U,\mathcal F)$ are 
defined as in \cite{DHdualcompl}, \S 2 and we set 
$$D^i(U,\mathcal F):=\ker [H^i(U,\mathcal F) \to\bigoplus_{v \not \in U} 
H^i(K_v,F)]={\rm Im} [H^i_c(U,\mathcal F) \to H^i(U,\mathcal F)],$$ 
where $F$ is the generic fibre of $\mathcal F$ over $K$.
For $v \in U$, we denote by $H^i_\textup{nr}(K_v,C_X)$
the image of $H^i(\calo_v,C_{\mathcal X})$ in $H^i(K_v,C_X)$ and if 
$S$ is a finite set of places that does not meet $U$, we set
$$\mathcal{P}^i_S(U,C_X) := \prod_{v \notin U, v \notin S} H^i(K_v,C_X) \times \prod_{v \in U} H^i_\textup{nr}(K_v,C_X)$$
(by convention, $\mathcal{P}^i_\emptyset(U,C_X)$ is denoted by 
$\mathcal{P}^i(U,C_X)$).

\begin{lem} \label{surjab}
	For every place $v \in U$, the abelianization map
	$\ab'_{{\mathcal X}_v} : \mathcal X(\calo_v) \to 
	H^0(\calo_v,C_{\mathcal X})$ is surjective. 
\end{lem}

\dem{} 

Using the same method as in \cite{DemEdinburgh},
Th.~2.18, it is sufficient (via the version of Lemma~\ref{lem compat ab ab'} over $\mathcal{O}_v$) to show that $H^1_{\ab}(\calo_v,\mathcal H)=0$
and $\ab^0 : \mathcal G(\calo_v) \to H^0_{\ab}(\calo_v,\mathcal G)$ is 
surjective. The nullity of $H^1_{\ab}(\calo_v,\mathcal H)$ follows 
by devissage from the fact that for an $\calo_v$-torus $\mathcal T$, we 
have $H^1(\calo_v,\mathcal T)=H^2(\calo_v,\mathcal T)=0$ (see for example 
\cite{HSz-Crelle}, proof of Th.~2.10). Finally 
$H^1(\calo_v, \mathcal G^{\sc}) \simeq H^1(\F_v, \widetilde{G^{\sc}})$
is trivial by Steinberg's Theorem (here $\widetilde{G^{\sc}}$ is the
reduction mod. $v$ of the reductive group scheme $\mathcal G^{\sc}$; it 
is a connected linear group scheme over the residue field $\F_v$ of 
the curve $E$ at $v$). This 
implies that the abelianization 
map $\mathcal G(\calo_v) \to H^0_{\ab}(\calo_v,\mathcal G)$ is  
surjective because by definition of the abelianization map, there is
an exact sequence 
$$\mathcal G(\calo_v) \to H^0_{\ab}(\calo_v,\mathcal G) \to H^1(\calo_v, \mathcal G^{\sc}).$$

\enddem{}

We need now to extend a few duality results of \cite{DHdualcompl} to
the three-term complex $C_X$~: 

\begin{prop} \label{dualcx}

a) The group $H^1_c(U, C_{\XX})$ is in the class $\mathcal{E}$ and 
$H^0(U,C_{\XX})$ is of finite type.

\smallskip

	b) The groups $D^0(U,C_{\XX})$ and
	$D^1(U, \widehat C_{\XX})$ are finite. 

\smallskip 

	c) The group $H^2_c(U,\widehat C_{\XX})$ is the dual of the 
	discrete group $H^0(U,C_{\XX})$ and 
	$H^1_c(U, C_{\XX})_\wedge \simeq H^1(U, \widehat{C}_{\XX})^D$.

\end{prop}

\dem{} a) Consider the exact triangle
        \begin{equation} \label{triangleC}
        \mu_H[2] \to C_{\XX} \to \mathcal C \to \mu_H[3] ,
        \end{equation}
        which implies that $H^1_c(U, C_{\XX})$ is an extension of
	$H^1_c(U,\mathcal C)$ by the finite (cf. \cite{DHAMM}, Th.~1.1) 
	group $H^3_c(U,\mu_H)$. Since
        $H^1_c(U,\mathcal C)$ is in $\mathcal{E}$
        (\cite{DHdualcompl}, Prop~3.13), so is $H^1_c(U, C_{\XX})$.

As $H^3(U,\mu_H)=0$ and $H^2(U,\mu_H)$ is finite 
(\cite{DHAMM}, Th.~1.1 and Cor.~4.9), we
also get that $H^0(U,C_{\XX})$ is an extension of $H^0(U,\mathcal C)$
(which is of finite type by \cite{DHdualcompl}, Prop~3.6. b) by a finite
group, hence it is also of finite type.

\smallskip

b) The finiteness of $D^0(U,C_{\XX})$ follows from that of
        $D^0(U,\mathcal C)$ (see \cite{DHdualcompl}, Lemma~4.15) and
        that of $H^2(U, \mu_H)$ (\cite{DHAMM}, Cor.~4.9)
        thanks to the commutative diagram with exact rows:

        $$
        \begin{CD}
        H^2(U,\mu_H) @>>> H^0(U,C_{\XX}) @>>> H^0(U,\mathcal C) \cr
        && @VVV @VVV \cr 
        && \bigoplus_{v \not \in U} H^0(K_v,C_X) @>>> \bigoplus_{v \not \in U}
        H^0(K_v,C)
        \end{CD}
        $$

Similarly $D^1(U, \widehat C_{\XX})$ contains $D^1(U,\widehat{\mathcal C})$
as a finite index subgroup thanks to exact triangle 
$$  \widehat {\mu_H} [-2] \to
\widehat {\mathcal C} \to \widehat {C_{\XX}} \to \widehat {\mu_H} [-1],$$
and $D^1(U,\widehat {\mathcal C})$ is finite by \cite{DHdualcompl}, Th.~5.2).

\smallskip

c) We first show that 
Artin-Verdier duality induces an isomorphism $H^1_c(U, C_{\XX})_\wedge \to H^1(U, \widehat{C}_{\XX})^D$. To prove this, we use a devissage given
        by the triangle (\ref{triangleC}).
        By \cite{DHdualcompl}, Theorem 4.11 b), the discrete torsion group
        $H^i(U,\widehat{\mathcal C})$ is dual to the profinite group
        $H^{2-i}(U, {\mathcal C})_{\wedge}$ for $i=1,2$. There is a
        commutative diagram
        $$
        \begin{CD}
        H^0_c(U,\mathcal C)_{\wedge} @>>> H^3_c(U,\mu_H) @>>>
        H^1_c(U,C_{\XX})_{\wedge} @>>> H^1_c(U,\mathcal C)_{\wedge} @>>> 0 \cr 
        @VVV @VVV @VVV @VVV \cr 
        H^2(U,\widehat{{\mathcal C}})^D @>>> H^0(U,\widehat \mu_H)^D @>>>
        H^1(U,\widehat{{\mathcal C}}_{\XX})^D @>>> H^1(U,\widehat{\mathcal C})^D
        @>>> 0 
        \end{CD}
        $$
        The second line is exact as the dual of an exact sequence
	of discrete torsion
        groups. The first line is exact as well thanks to Lemma~3.12 of
        \cite{DHdualcompl} because $H^3_c(U,\mu_H)$ is finite
        and $H^1_c(U,C_{\XX})$ belongs to the class $\mathcal E$ (hence
        it injects into its completion).
        Now Artin-Mazur-Milne duality for $\mu_H$
	(see \cite{DHAMM}, Theorem 1.1) and
        for $\mathcal C$ (see \cite{DHdualcompl}, Theorem 4.11 b)) yield that
        the first, second, and fourth vertical maps are isomorphisms.
        The five lemma implies that the third one is also an isomorphism, as
        required.

\smallskip

The argument to prove the isomorphism $H^1_c(U, C_{\XX})_\wedge \simeq H^1(U, \widehat{C}_{\XX})^D$, is similar, using the commutative diagram with 
exact rows:
$$
        \begin{CD}
0 @>>>	H^2_c(U,\widehat{\mathcal C}) @>>> H^2_c(U,\widehat{\mathcal C_{\XX}})
	@>>> H^1_c(U,\widehat \mu_H) @>>> H^3_c(U,\widehat {\mathcal C})\cr
		&&      @VVV @VVV @VVV @VVV \cr
		0 @>>> H^0(U,\mathcal C)^D @>>> H^0(U,\mathcal C_{\XX})^D @>>>
		H^2(U,\mu_H)^D @>>> H^{-1}(U,{\mathcal C})^D 
        \end{CD}
        $$

	Indeed the first and fourth vertical maps are isomorphisms by 
	\cite{DHdualcompl}, Th. 4.9 a), 
	and so is the third by \cite{DHAMM}, Theorem 1.1.

\enddem

\begin{prop} \label{ptsuit}
	
a) There is an exact sequence 
$$H^0(K, C_X) \to H^0(\A_K, C_X) \to H^1(K, \widehat{C}_X)^D .$$

\smallskip

	b) There is an exact sequence 
$$H^1(K,\widehat C_X) \to H^1(\A_K, \widehat C_X) \to 
H^0(K,C_X)^D.$$

\smallskip

c) There is a 
local duality isomorphism
$$H^0(K_v, C_X)_\wedge \xrightarrow{\sim} H^1(K_v, \widehat{C}_X)^D,$$
and the orthogonal of $H^1_{\nr}(K_v, \widehat{C}_X)$ in the local duality 
is the image of $H^0_{\nr}(K_v, C_X)$ in $H^0(K_v, C_X)_\wedge$.

\end{prop}

\dem{} a) We extend the proof of the exactness of the second line
of Poitou-Tate exact sequence in Theorem 5.7 in \cite{DHdualcompl}.
Following the same method, it is sufficient to
extend the first part of
Lemma 5.6 a) in \cite{DHdualcompl} to the complex $C_X$. Namely,
it remains to show the exactness of 

	$$H^0(U,C_{\XX}) \to \mathcal{P}^0(U,C_X) \to
        \left(H^1(K, \widehat{C}_X) / \Sha^1(K, \widehat{C}_X)\right)^D .$$
The proof is exactly the same as for $C$, applying 
Lemma 2.2 in \cite{DHdualcompl} to the complex $C_X$ and using the three
following facts (proven in Proposition~\ref{dualcx}): for every 
non empty Zariski open subset $V \subset U$, 
the group $H^1_c(V, C_{\XX})$ is in the class $\mathcal{E}$,  
the group $D^0(V,C_{\XX})$ is finite, and Artin-Verdier duality 
induces an isomorphism
$H^1_c(V, C_{\XX})_\wedge \simeq H^1(V, \widehat{C}_{\XX})^D$.

\smallskip

b) Similarly, it is sufficient to extend the second part of Lemma~5.6 b) 
in \cite{DHdualcompl}, that is to show the exactness of 
$$H^1(U, \widehat C_{\XX}) \to \mathcal{P}^1(U, \widehat C_X)
\to H^0(K,C_X)^D.$$
Applying again Lemma 2.2 of \cite{DHdualcompl} to $\widehat C_X$, one just has 
to check the following properties~: 

-for every non empty Zariski open subset $V \subset U$, 
we have $$H^2_c(V,\widehat C_{\XX}) \simeq (H^0(V,C_{\XX})_{\wedge})^D = 
(H^0(V,C_{\XX})^D.$$ This holds thanks to Proposition~\ref{dualcx} (the 
finite type group $H^0(V,C_{\XX})$ has same dual as $H^0(V,C_{\XX})_{\wedge}$).

-the group $D^1(U,\widehat C_{\XX})$ is finite, which is also proven in 
Proposition~\ref{dualcx}.

\smallskip

c) Using the exact triangle (\ref{triangleC}) and
the vanishing of $H^3(K_v,\mu_H)$ (\cite{MilADT}, Prop. III.6.4),
there is a commutative diagram with exact rows (the completed
first line remains exact because $H^2(K_v,\mu_H)$ is finite by \cite{MilADT},
Ex. III.6.7, and the second line is obtained by dualizing an exact sequence
of discrete torsion groups):
$$
\begin{CD}
H^{-1}(K_v,C)_{\wedge }@>>> H^2(K_v,\mu_H)@>>>
H^0(K_v,C_X)_{\wedge} @>>> H^0(K_v,C)_{\wedge}  @>>> 0 \cr
@VVV @VVV @VVV @VVV \cr
H^2(K_v,\widehat{C})^D @>>> H^0(K_v,\widehat \mu_H)^D @>>> 
H^1(K_v, \widehat{C}_X)^D @>>> H^1(K_v, \widehat{C})^D @>>> 0
\end{CD}
$$
Since the first, second, and fourth vertical map are isomorphisms by
\cite{MilADT}, Th.~III.6.10 and \cite{Demtores}, Th.~3.1., so is the third
vertical map.

\smallskip

It remains to show that the map $H^0(K_v,C_X)_{\wedge}/H^0_{\nr}(K_v,C_X) 
\to H^1(\calo_v,\widehat{C}_{\XX})^D$ is injective. There is a commutative 
diagram with exact lines:

$$
\begin{CD}
	&&&& 
H^0(\calo_v,C_{\XX}) @>>> H^0(\calo_v,\C)  @>>> 0 \cr
	&&&&  @VVV @VVV \cr
H^{-1}(K_v,C)_{\wedge }@>>> H^2(K_v,\mu_H)@>>>
H^0(K_v,C_X)_{\wedge} @>>> H^0(K_v,C)_{\wedge}  @>>> 0 \cr
@VVV @VVV @VVV @VVV \cr
H^2(\calo_v,\widehat{\C})^D @>>> H^0(\calo_v,\widehat \mu_H)^D @>>>
	H^1(\calo_v, \widehat{C}_{\XX})^D @>>> H^1(\calo_v, \widehat{\C})^D @>>> 0
\end{CD}
$$

The fourth column is exact and the map 
$H^{-1}(K_v,C)_{\wedge } \to H^2(\calo_v,\widehat{\C})^D$ is surjective 
(\cite{Demtores}, Th. 3.1. and 3.3). Since $\widehat \mu_H$ is a finite 
group scheme, we have 
$H^0(\calo_v,\widehat \mu_H)=H^0(K_v,\widehat \mu_H)$, hence 
the map $H^2(K_v,\mu_H) \to H^0(\calo_v,\widehat \mu_H)^D$ is injective 
by \cite{MilADT}, Th.~III.6.10. The required result follows by diagram 
chasing.

\enddem

\begin{prop} \label{dualstep}
	
	There is an exact sequence
	\begin{equation} \label{short exact seq complex SA}
	H^0(K, C_X) \to H^0(\A_K, C_X) \to \left(H^1(K, \widehat{C}_X) / \Sha^1(K, \widehat{C}_X)\right)^D \to {H^{-1}(K, \widehat{C}_X)^\star}_\wedge / H^{-1}(K, \widehat{C}_X)^\star \to 0, 
	\end{equation}
	where the last non trivial map is defined via the natural map 
	$$\left(H^1(K, \widehat{C}_X) / \Sha^1(K, \widehat{C}_X)\right)^D \to 
	\left(H^1(K, \widehat{C}) / \Sha^1(K, \widehat{C})\right)^D $$ 
	and the map $$\left(H^1(K, \widehat{C}) / \Sha^1(K, \widehat{C})\right)^D
	\stackrel{\partial_C}{\to} 
	{H^{-1}(K, \widehat{C})^\star}_\wedge / H^{-1}(K, \widehat{C})^\star 
	\simeq {H^{-1}(K, \widehat{C_X})^\star}_\wedge / H^{-1}(K, \widehat{C_X})^\star.$$

\end{prop}

\dem{} We consider diagram \eqref{diag hom spaces devissage}
comparing $C_X$ and $C$. Applying cohomology, we get a commutative diagram
\begin{equation} \label{devissage complex 3-terms}
\xymatrix{
	H^2(\A_K, \mu_H) \ar[r] \ar[d] & \left(H^0(K, \widehat{\mu_H})\right)^D \ar[r] \ar[d] & 0 &  \\
	H^0(\A_K, C_X) \ar[r] \ar[d] & \left(H^1(K, \widehat{C}_X) / \Sha^1(K, \widehat{C}_X)\right)^D \ar[r] \ar[d] & {H^{-1}(K, \widehat{C}_X)^\star}_\wedge / H^{-1}(K, \widehat{C}_X)^\star \ar[r] \ar[d]^{\simeq} & 0 \\
	H^0(\A_K, C) \ar[r] \ar[d] & \left(H^1(K, \widehat{C}) / \Sha^1(K, \widehat{C})\right)^D \ar[r] \ar[d] & {H^{-1}(K, \widehat{C})^\star}_\wedge / H^{-1}(K, \widehat{C})^\star \ar[r] & 0 \\
	0 & 0 & 
}
\end{equation}
Local duality for $\mu_H$ (cf. \cite{Ces}, Proposition 4.10 b) for instance) and 
Proposition~\ref{dualstep} imply that the first row and the last one are exact. The exact triangles of \eqref{diag hom spaces devissage} implies that the first 
column is exact.

The second column is exact, being the dual of the obvious exact sequence
\[
0 \to H^1(K, \widehat{C}) / \Sha^1(K, \widehat{C}) \to H^1(K, \widehat{C}_X) / \Sha^1(K, \widehat{C}_X) \to H^0(K, \widehat{\mu_H}) \, .
\]

Now an easy diagram chase implies the exactness of the second line of \eqref{devissage complex 3-terms}.
The exactness of \eqref{short exact seq complex SA} then follows from 
Proposition~\ref{ptsuit}, a).

\enddem

\begin{prop} \label{flatcx}

	a) The following sequence:
	\begin{equation} \label{exact seq compl U 3}
	H^0(K,C_X)_\wedge \to H^0(\A_K,C_X)_\wedge \to \left(H^1(K, \widehat{C}_X) / \Sha^1(K, \widehat{C}_X)\right)^D \to 0
	\end{equation}

	is exact.
	
	\smallskip
	
	b) Let $S$ be a finite set of places that does not meet $U$. Set 
	$$H^1_S(K, \widehat{C}_X) =\ker [H^1(K, \widehat{C}_X) 
	\to \bigoplus_{v \in S} H^1(K_v, \widehat{C}_X)].$$
	Then the sequence
	\begin{equation} \label{exact seq compl U S 3}
	H^0(K,C_X)_\wedge \to H^0(\A_K ^S,C_X)_\wedge \to \left(H^1_S(K, \widehat{C}_X)/\Sha^1(K, \widehat{C}_X)\right)^D \to 0
	\end{equation}
	is exact.

	\smallskip

	c) There is an exact sequence 
	$$H^0(U,C_{\XX})_{\wedge} \to {\bf P}^0_S(U,C_{\XX})_{\wedge} \to 
	\left(H^1_S(K, \widehat{C}_X)/\Sha^1(K, \widehat{C}_X)\right)^D.$$
	
\end{prop}

\dem{} a) Since $\left(H^1(K, \widehat{C}_X) /
\Sha^1(K, \widehat{C}_X)\right)^D$ is profinite and
${H^{-1}(K, \widehat{C}_X)^\star}_\wedge / H^{-1}(K, \widehat{C}_X)^\star$
is uniquely divisible, the surjectivity of the map 
$H^0(\A_K,C_X)_\wedge \to \left(H^1(K, \widehat{C}_X) / \Sha^1(K, \widehat{C}_X)\right)^D$ follows immediately from Proposition~\ref{dualstep}. Applying 
Proposition~\ref{ptsuit} b) 
and taking into account that $H^1(K,C_X)$ is torsion, we get
an exact sequence
$$0 \to \Sha^1(K,\widehat C_X) \to 
H^1(K,\widehat C_X) \to H^1(\A_K, \widehat C_X)_{\tors} \to
(H^0(K,C_X)_{\wedge})^D.$$
Indeed $(H^0(K,C_X)_{\wedge})^D=(H^0(K,C_X)^D)_{\tors}$ (cf.
\cite{DHdualcompl}, Remark~5.9). Dualizing this exact sequence  now 
yields the result, thanks to Proposition~\ref{ptsuit}, c).

\smallskip

b) Dualizing the exact sequence of discrete torsion groups 
\[
0 \to H^1_S(K, \widehat{C}_X) \to H^1(K, \widehat{C}_X) \to \bigoplus_{v \in S} H^1(K_v, \widehat{C}_X) \, .
\]
and taking into account the local duality isomorphism
$H^0(K_v, C_X)_\wedge \xrightarrow{\sim} H^1(K_v, \widehat{C}_X)^D$, one gets a commutative diagram
\[
\xymatrix{
	& \bigoplus_{v \in S} H^0(K_v, C_X)_\wedge \ar[r]^= \ar[d]^i & \bigoplus_{v \in S} H^0(K_v, C_X)_\wedge \ar[d] & \\
	H^0(K,C_X)_\wedge \ar[r] \ar[d]^= & H^0(\A_K,C_X)_\wedge \ar[r] \ar[d]^p & \left(H^1(K, \widehat{C}_X)/\Sha^1(K, \widehat{C}_X)\right)^D \ar[d] \ar[r] & 0 \\
	H^0(K,C_X)_\wedge \ar[r] & H^0(\A_K^S,C_X)_\wedge \ar[r] \ar[d] & \left(H^1_S(K, \widehat{C}_X)/\Sha^1(K, \widehat{C}_X)\right)^D \ar[d] & \\
	& 0 & 0 \, ,
}
\]
where the top line (by a) and the columns are exact (the left one 
because $H^0(K_v,C_X)$ is in the class $\mathcal E$ by
\cite{DHdualcompl}, Prop~3.13 and exact triangle (\ref{triangleC});
so the sequence remains 
exact after completion thanks to loc. cit., Lemma~3.12). 
A simple diagram chase implies that the bottom line is exact, which 
proves b).

\smallskip

c) For $v \in U$ and $i \in \Z$, set
$$H^i_r(K_v,\widehat C_X)=H^i(K_v,\widehat C_X)/H^i_{\nr}(K_v,\widehat C_X)$$
(and similarly for $C_X$).
Consider the commutative diagram: 
\[
\xymatrix{
	H^1_S(K,\widehat C_X) \ar[r] &  \bigoplus_{v \in U} H^1_r(K_v,\widehat C_X)
	\oplus \bigoplus_{v \not \in (S \cup U) }
	H^1(K_v,\widehat C_X) \ar[r] & H^0(U,C_{\XX})^D \\
	H^1_S(K,\widehat C_X) \ar[r] \ar[u]^= & \prod'_{v \not \in S} H^1(K_v, \widehat C_X) \ar[r] \ar[u]^p & 
	H^0(K,C_X)^D \ar[u] \\
	& \prod_{v \in U} H^1_{\nr}(K_v,\widehat C_X) \ar[u] \ar[r]^= & \prod_{v \in U} 
	H^1_{\nr}(K_v,\widehat C_X) \ar[u] \, .
}
\]
The middle column is obviously exact. The right column is exact as well,
because by Proposition~\ref{ptsuit} c),
it is the dual of the sequence of discrete groups
$$H^0(U,C_{\XX}) \to H^0(K,C_X) \to \bigoplus_{v \in U} H^0_r(K_v,C_X),$$
which is exact by \cite{DHdualcompl}, Prop.~2.1. It follows immediately 
from Proposition~\ref{ptsuit} b) that the second line of the diagram is 
exact. Since the map $p$ is obviously surjective, a diagram chase shows that 
the first line of the diagram is exact as well. Dualizing it (and observing 
that $H^0(U,C_{\XX})$ is an abelian group of finite type by
Proposition~\ref{dualcx}, a) yields
(thanks to Proposition~\ref{ptsuit}, c) the exactness of 
$$H^0(U,C_{\XX})_{\wedge} \to {\bf P}^0_S(U,C_{\XX})_{\wedge} \to
        H^1_S(K, \widehat{C}_X)^D.$$
Since the canonical map ${\bf P}^0_S(U,C_{\XX})_{\wedge} \to 
H^1_S(K, \widehat{C}_X)^D$ factors through $(H^1_S(K, \widehat{C}_X)/\Sha^1(K,\widehat{C}_X))^D$, the result is proven.

\enddem

Recall that for a finite (possibly empty) 
set of places $S$ of $K$, we have the Brauer-Manin
pairing 
$${\rm BM} : X(A_K ^S) \times \br X \to \Q/\Z , \quad 
((P_v)_{v \not \in S} , \alpha) \mapsto \sum_{v \not \in S} \alpha(P_v).$$
By global class field theory, elements of $X(K) \subset X(A_K ^S)$
are orthogonal to 
$\Br_S(X)$ for this pairing; in particular when $S=\emptyset$, elements
of $X(K)$ are orthogonal to $\Br_e X$ (hence to $\br X$).
By continuity of the pairing,
the same holds for the closure $\overline{{X(K)}}^S$ of 
$X(K)$ in $X(A_K ^S)$ for the strong topology. The following theorem 
gives various converse statements. 

\begin{theo} \label{mainstrong}
	Let $X=G/H$ be a homogeneous space of a reductive group $G$ with 
	$H$ reductive. Set $U(X):=K[X]^*/K^*$. 
	Assume that $G^\sc$ satisfies strong approximation
	outside $S_0$ (finite set of places).
	\begin{enumerate}
		\item There is a natural exact sequence of pointed topological spaces:
		\[
		1 \to \rho\left(\overline{G^\sc(K) \cdot G^\sc(K_{S_0})}\right) \cdot X(K) \to X(\A_K) \xrightarrow{BM} \left(\Br_{1,e}(X,G)/\Be(X)\right)^D \xrightarrow{\partial} {U(X)^\star}_\wedge / U(X)^\star \to 1 \, .
		\]
		In particular, $X(\A_K)^{\Br X} = \rho\left(\overline{G^\sc(K) \cdot G^\sc(K_{S_0})}\right) \cdot X(K)$ and $X(\A_K^{S_0})^{\Br X} \subset
		\overline{X(K)}^{S_0}$.
		\item If $S$ is a non empty finite set of places, there is a natural exact sequence of pointed topological spaces:
		\[
		1 \to \overline{\rho(G^\sc(K_{S_0})) \cdot X(K)}^S \to X(\A_K^S) \xrightarrow{BM} \left(\Br_S(X,G)/\Be(X)\right)^D \to 1 \, .
		\]
		In particular, $X(\A_K^S)^{\Br_S X} =
		\overline{\rho(G^\sc(K_{S_0})) \cdot X(K)}^S$.
		
		If $S_0 \subset S$, we get an exact sequence
		\[
		1 \to \overline{X(K)}^S \to X(\A_K^S) \xrightarrow{BM} \left(\Br_S(X,G)/\Be(X)\right)^D \to 1 \, .
		\]
	\end{enumerate}
\end{theo}

\dem{} As earlier, one can assume that the group $G$ is quasi-trivial, 
up to replacing $G$ by a flasque resolution
$$1 \to S \to G' \to G \to 1.$$ 
Indeed $\pic \ov S=0$ (since $S$ is a torus), hence $\br \ov G \hookrightarrow
\br \ov G'$, which implies $\br_1(X,G)=\br_1(X,G')$ by \cite{San},
Prop.~6.10. 
Throughout the proof, we use Theorem~\ref{theo compatible strong approx} to translate results concerning complexes of tori and cup-products to results concerning homogeneous spaces and Brauer-Manin obstructions.

\begin{enumerate}
	
	\item By \cite{DemJOA}, Th. 4.14, the group $U(X)$ is isomorphic to
	$H^{-1}(K,\widehat{C_X})$. 
	By Theorems~\ref{theo compatible strong approx} and 
	\ref{dualstep}, there is a commutative 
	diagram with exact second row: 
	{\small 
		$$
		\begin{CD}
		X(K) @>>> X(\A_K) @>{\rm BM}>> \left(\Br_{1,e}(X,G)/\Be(X)\right)^D @>>>
		{U(X)^\star}_\wedge / U(X)^\star @>>> 1 \cr 
		@VV{\ab ' _K }V @VV{\ab ' _{\A_K}}V @V{\simeq}V{(\Phi_X)^D}V
		@VV{\simeq}V \cr 
		H^0(K,C_X) @>>> H^0(\A_K,C_X) @>>> H^1(K,\widehat{C_X})^D @>>> 
		{H^{-1}(K, \widehat{C_X})^\star}_\wedge / H^{-1}(K, \widehat{C_X})^\star 
		@>>> 1
		\end{CD}
		$$
	}

		By Proposition~\ref{orbitesprop} and 
		Lemma~\ref{surjab},
		the maps $\ab ' _K$ and $\ab '_{\A_K}$ are surjective.
	By diagram chasing, the sequence of three last non-trivial terms on the 
	first line is also exact. Besides, every element $x$ of $X(K) \subset X(\A_K)$ 
	satisfies ${\rm BM}(x)=0$ by class field theory, and the same holds 
	for an element of $X(\A_K)$ of the form $(g_v.x)$ 
	with $(g_v) \in \rho(G^{\sc}(\A_K))$ thanks to the commutativity of 
	the diagram and the property $\ab '(g_v.x)=\ab '(x)$
		for every place $v$ of $K$ (Proposition~\ref{prop ab action}).
	It remains to show conversely that 
	every $(P_v) \in X(\A_K)$ such that ${\rm BM}((P_v))=0$
	comes from
	$\rho\left(\overline{G^\sc(K) \cdot G^\sc(K_{S_0})}\right) \cdot X(K)$ 
	The diagram and the surjectivity of $\ab ' _K$ imply that there exists 
	$x \in X(K)$ such that $$\ab '_{\A_K}(x)=\ab '_{\A_K}(P_v) 
	\in H^0(K_v,C_X)$$ for every place $v$.
		By Proposition~\ref{orbitesprop},
	there exists for each $v$ an element $g_v \in \rho(G^{\sc}(K_v))$
	such that $P_v=g_v.x$, and we can assume 
	that $g_v \in {\mathcal G}^{\sc}(\calo_v)$ for $v \not \in S$,
	where $S \supset S_0$ is some finite set of places. Since 
	$G^{\sc}$ satisfies strong approximation outside $S_0$,
	we finally obtain that $(P_v)$ belongs to 
	$\rho\left(\overline{G^\sc(K) \cdot G^\sc(K_{S_0})}\right) \cdot X(K)$ as 
	required. The two other assertions in 1. follow immediately.

	\item We follow the proof of \cite{DemEdinburgh}, Theorem 6.1.
	
	Up to shrinking $U$, we can assume that it does not meet $(S_0 \cup S)$.
	We consider the following commutative diagram:
	\begin{changemargin}{-1cm}{1cm}
		\begin{equation} \label{big diag}
		\xymatrix{
			& & & \mathcal{P}_S(U,C_H)_\wedge \ar[rr] \ar'[d][dd] & & H^1_S(K, \widehat{C}_X)^D \ar[dd] \\
			& & \mathcal{P}_S(U, H) \ar[ur] \ar[rr] \ar[dd] & & \left(\Br_S(H) \right)^D \ar[ru] \ar[dd] & \\
			& H^0(U,C_G)_\wedge \ar'[d][dd] \ar'[r][rr] & & \mathcal{P}^0_S(U,C_G)_\wedge \ar'[r][rr] \ar'[d][dd] & & H^1_S(K, \widehat{C}_G)^D \ar[dd] \\
			G(U) \ar[ru] \ar[rr] \ar[dd] & & \mathcal{P}^0_S(U,G) \ar[rr] \ar[dd] \ar[ru] & & \left(\Br_S(G) \right)^D \ar[ru] \ar[dd] & \\
			& H^0(U, C_X)_\wedge \ar'[r][rr] \ar'[d][dd] & & \mathcal{P}^0_S(U, C_X)_\wedge \ar'[r][rr] \ar'[d][dd] & & H^1_S(K, \widehat{C}_X)^D \\
			X(U) \ar[ru] \ar[rr] \ar[dd] & & \mathcal{P}^0_S(U,X) \ar[ru] \ar[rr] \ar[dd] & & \left(\Br_S(X,G)\right)^D \ar[ru] & \\
			& H^1(U, C_H) \ar'[r][rr] & & \mathcal{P}^0_S(U, C_H) & & \\
			H^1(U,H) \ar[d] \ar[rr] \ar[ru]^\sim & & \mathcal{P}^1_S(U,H) \ar[ru]^\sim & & & \\
			1 & & & & & \, .}
		\end{equation}
	\end{changemargin}
	The definition of the maps $H^0(U,C_X)_\wedge \to H^1(U,C_H)$ and $\mathcal{P}^0_S(U,C_X)_\wedge \to \mathcal{P}^1_S(U,C_H)$ is a consequence of Lemma 6.5 in \cite{DemEdinburgh} (the proof works the same in the function field context). Using this lemma, we see that the columns in this diagram are exact. The surjectivity of $X(U) \to H^1(U,H)$ is a consequence of the vanishing of $H^1(U,G)$ for $U$ sufficiently small (Hilbert 90, together with Nisnevich Theorem as proven in \cite{Gil-torseurs-affine}, Theorem 5.1).
	
	In addition, the fifth line of this diagram is exact by
	Proposition~\ref{flatcx}, c), and the maps $H^1(U,H) \to H^1(U, C_H)$
	and $\mathcal{P}^1_S(U,H) \to \mathcal{P}^1_S(U,C_H)$ are bijections (see \cite{GA}  Theorem 5.5 and example 5.4 (iii), and \cite{Gil-torseurs-affine}, Theorem 5.1).
	
	Now, the proof of the second point of the Theorem is a diagram chase
	in diagram \eqref{big diag}, following exactly the argument in
	\cite{DemEdinburgh}, Theorem 6.1. In particular, one uses the fact that
	$\mathcal{P}^0_S(U,C_H)_\wedge$ and $\mathcal{P}^0_S(U,C_H)$
	(hence $\mathcal{P}^0_S(H)$) have the same image inside
	$H^1_S(K, \widehat{C}_H)^D$ which is a consequence of the following 
	proposition~: 
	
	\begin{prop} \label{compactquot} 
		
		Let $\mathcal D=[\T_1 \to \T_2]$ be a complex 
		(with $\T_2$ in degree $0$) of flat,
		separated and finite type commutative group schemes over 
		the affine curve $E-S$, such that the restriction of 
		$\T_1$ and $\T_2$ to $U$ are tori. Let $D=[T_1 \to T_2]$
		be the generic fibre of $\mathcal D$.
		
		\smallskip
		
		a) The group $H^0(U,\D) \backslash \mathcal{P}^0_S(U,D) /
		\prod_{v \notin S} H^0(\mathcal{O}_v,\D)$ is finite.
		
		\smallskip
		
		b) The group $\mathcal{P}^0_S(U,\D) / \overline{H^0(U,\D)}$ is compact.
		
	\end{prop}
	
	\dem{} a) Let $S_1$ be the finite set of places of $K$ that are neither 
	in $U$ nor in $S$. Recall (\cite{MilADT}, Cor.~2.3) 
	that the group $H^1(K_v,T)$ is finite for 
	all places $v \in \Omega_K$.
	If $\T_1=0$, the result follows immediately from the finiteness of 
	the class group of a torus, as proven in \cite{Conrad-Compositio}, Example 1.3.2 (recall that $S \neq \emptyset$). In the general case, we have $H^1(\calo_v,\T_1)=0$ for $v \in U$, whence
	an exact sequence 
	$$\mathcal{P}^0_S(U,\T_2) \to \mathcal{P}^0_S(U,D) \to
	\prod_{v \in S_1} H^1(K_v,T_1),$$
	which shows in particular 
	that the image of $\mathcal{P}^0_S(U,\T_2)/\prod_{v \notin S} H^0(\mathcal{O}_v,\T_2)$ is of finite index in 
	$\mathcal{P}^0_S(U,\D)/\prod_{v \notin S} H^0(\mathcal{O}_v,\D)$. Thus one 
	reduces to the already known case $\T_1=0$.
	
	\smallskip 
	
	b) follows from a) and the compactness of
	$\prod_{v \notin S} H^0(\mathcal{O}_v,\D)$.
	\enddem
	
	To finish the proof of the Theorem, one deduces the surjectivity of the map
	$$X(\A_K^S) \to \left(\Br_S(X,G) / \Be(X)\right)^D$$ from exact sequence \eqref{exact seq compl U S 3} and from the surjectivity of $X(\A_K^S) \to H^0(\A_K^S,C_X)$ (which follows from Lemma~\ref{surjab} and
	Proposition~\ref{orbitesprop}).

\end{enumerate}
\enddem

\begin{rema} {\rm
	\begin{enumerate}
		\item In the special case when $X=G$ is a reductive group, we get exact
sequences of topological groups: 

	$$ 1 \to \rho\left(\overline{G^\sc(K) \cdot G^\sc(K_{S_0})}\right) \cdot G(K) \to G(\A_K) \xrightarrow{BM} \left(\Br_{1,e} \, G/\Be(G)\right)^D \xrightarrow{\partial} {\widehat G(K)^\star}_\wedge / \widehat G(K)^\star \to 1 \, .$$

               $$ 
		1 \to \overline{\rho(G^\sc(K_{S_0})) \cdot G(K)}^S \to G(\A_K^S) \xrightarrow{BM} \left(\Br_{1,S} \, G /\Be(G)\right)^D \to 1 \, .
               $$ 
The fact that $BM$ is a group morphism in this situation follows from the 
	same argument as in Corollary~\ref{groupcase}.
	\item In addition, using Theorem 5.10 in \cite{DHdualcompl}, it is straightforward to continue the previous exact sequence as follows, for any non-empty finite set of places $S$:
\[
\xymatrix{
1 \ar[r] & \overline{\rho(G^\sc(K_{S_0})) \cdot G(K)}^S \ar[r] & G(\A_K^S) \ar[r]^{BM} & (\Br_{1,S} \, G)^D \ar[d] \\
& (\pic G)^D & \bigoplus_v H^1(K_v,G) \ar[l] & H^1(K,G) \ar[l]  \, ,
}
\]		
where the pairing $H^1(K_v,G) \times \pic(G)$ can be defined via the natural isomorphism between $\pic(G)$ and the group 
			$\textup{Ext}^1_K(G, \Gm)$ of central extensions of $G$ by $\Gm$ (see for instance Corollary 5.7 and Proposition 8.2 in \cite{CTfl}). One can even extend this exact sequence to a 9-term Poitou-Tate exact sequence involving non-abelian $H^2$, using Theorem 5.10 in \cite{DHdualcompl}, following the methods of \cite{DemAF}, Theorem 5.1. See also \cite{BorAMS}, Th~5.16 for a
similar statement over a number field.
\end{enumerate}}
\end{rema}

\bigskip

\noindent 
Sorbonne Universit\'e and Universit\'e de Paris, CNRS, IMJ-PRG, F-75006 Paris, France. 

\noindent
cyril.demarche@imj-prg.fr

\bigskip

\noindent 
Universit\'e Paris-Saclay, CNRS, Laboratoire de math\'ematiques d'Orsay, 91405, Orsay, France.

\noindent 
david.harari@math.u-psud.fr


\begin{thebibliography}{}

\bibitem[Bor1]{BorDuke} M. Borovoi: {\it Abelianization of the second nonabelian Galois cohomology}, Duke Math. J., {\bf 72}, 217--239 (1993).

\bibitem[Bor2]{BorAMS} M. Borovoi: {\it Abelian {G}alois cohomology of reductive groups}, Mem. Amer. Math. Soc. {\bf 132}, no. 626 (1998).

\bibitem[Bor3]{Borcrelle} M. Borovoi: {\it The Brauer-Manin obstructions for homogeneous spaces with connected or abelian stabilizer},
	J. Reine Angew. Math. {\bf 473}, 181–-194 (1996).

\bibitem[Bor4]{borann} M. Borovoi: {\it A cohomological obstruction to
the Hasse principle for homogeneous spaces}, Math. Ann. {\bf 314},
491--504 (1999).

\bibitem[BD]{BDcom} M. Borovoi, C. Demarche~: {\it Manin obstruction to strong approximation for homogeneous spaces},
	Comment. Math. Helv. {\bf 88}, no. 1, 1–-54 (2013).

\bibitem[BDH]{BDH} M. Borovoi, C. Demarche, D. Harari: {\it Complexes de groupes de type multiplicatif et groupe de Brauer non ramifi\'e des espaces homog\`enes}, Ann. Sci. \'Ec. Norm. Sup\'er. {\bf 46}, No. 4, 651--692 (2013).

\bibitem[BvH]{BvH} M. Borovoi, J. van Hamel: {\it Extended Picard complexes and linear algebraic groups}, J. reine angew. Math. {\bf 627}, 53--82 (2009).

\bibitem[BvH2]{BvH2} M. Borovoi, J. van Hamel: {\it Extended equivariant Picard complexes and homogeneous spaces}, , Transform. Groups {\bf 17}, 51--86 (2012).

\bibitem[BT]{BT} F. Bruhat,  J. Tits: {\it Groupes  alg\'ebriques  sur  un  corps  local. Chapitre  III. Compl\'ements et applications \`a la cohomologie galoisienne.} J. Fac. Sci. Univ. Tokyo Sect. IA Math. {\bf 34}, no. 3, 671--698 (1987). 

\bibitem[CT]{CTfl} J.-L. Colliot-Th\'el\`ene: {\it R\'esolutions flasques des groupes lin\'eaires connexes}, J. f\"ur die reine und angewandte Mathematik {\bf 618}, 77--133 (2008). 

\bibitem[CTS]{CTSdesc} J.-L. Colliot-Th\'el\`ene, J.-J. Sansuc: {\it La descente sur les vari\'et\'es rationnelles. II}, 
Duke Math. J. {\bf 54} no. 2, 375--492 (1987).

\bibitem[CTX]{CTX} J.-L. Colliot-Th\'el\`ene, F. Xu: {\it 
Brauer-Manin obstruction for integral points of homogeneous spaces and representation by integral quadratic forms}, Compos. Math. {\bf 145} no. 2, 309–-363 (2009).

\bibitem[Con]{Conrad-Compositio} B. Conrad: {\it Finiteness theorems for algebraic groups over function fields}, Compos. Math. {\bf 148}, No. 2, 555--639 (2012).

\bibitem[DeB]{DeB} S. DeBacker: {\it Parameterizing conjugacy classes of maximal unramified tori via Bruhat-Tits theory}, Michigan Math. J. 54, 157--178 (2006).

\bibitem[Ces]{Ces} K. \v{C}esnavi\v{c}ius: {\it Poitou–-Tate without restrictions on the order}, Mathematical Research Letters {\bf 22}, no. 6, 1621--1666 (2015).

\bibitem[Dem1]{Demtores} C. Demarche: {\it Suites de Poitou-Tate pour les complexes de tores \`a deux termes}, Int. Math. Res. Notices, Vol. {\bf 2011} (1), 135--174 (2011).

\bibitem[Dem2]{DemAF} C. Demarche: {\it Le d\'efaut d'approximation forte dans les groupes lin\'eaires connexes}, Proc. London Math. Soc. {\bf 102} (3), 563--597 (2011).

\bibitem[Dem3]{DemJOA} C. Demarche: {\it Une formule pour le groupe de Brauer alg\'ebrique d'un torseur}, Journal of Algebra {\bf 347} (1), 96--132 (2011).

\bibitem[Dem4]{DemEdinburgh} C. Demarche: {\it Ab\'elianisation des espaces homog\`enes et applications arithm\'etiques}, in "Torsors, \'etale homotopy and applications to rational points". Cambridge University Press, London Mathematical Society Lecture Note Series 405, 138--209 (2013). 

\bibitem[DH1]{DHAMM} C. Demarche, D. Harari: {\it Artin-Mazur-Milne duality for fppf cohomology}, Algebra Number Theory {\bf 13}, No. 10, 2323--2357 (2019). 

\bibitem[DH2]{DHdualcompl} C. Demarche, D. Harari: {\it Duality for complexes of tori over a global field of positive characteristic},
Journal de l'École Polytechnique (Mathématiques) {\bf 7},
831--870 (2020). 

\bibitem[Demeio]{Demeio} J. Demeio~: {\it The étale Brauer-Manin obstruction 
	to strong approximation on homogeneous spaces}, Preprint 2021.

\bibitem[FSS]{FSS} Y. Flicker, C. Scheiderer, R. Sujatha~: 
	{Grothendieck's theorem on non-abelian $H^2$ and local-global principles},
		J. Amer. Math. Soc. {\bf 11}, no. 3, 731–-750 (1998).

\bibitem[Gil]{Gil-torseurs-affine} P. Gille: {\it Torseurs sur la droite affine}, Transform. Groups {\bf 7}, No. 3, 231--245 (2002).

\bibitem[GA]{GA} C. Gonz\'alez-Avil\'es: {\it Quasi-abelian crossed modules and nonabelian cohomology}, J. Algebra {\bf 369}, 235--255 (2012).

\bibitem[HSk]{HS-op-desc} D. Harari, A. Skorobogatov: {\it Descent theory for open varieties}, in "Torsors, \'etale homotopy and applications to rational points", London Mathematical Society Lecture Note Series 405, 250--279 (2013).

\bibitem[HSz1]{HSz-Crelle} D. Harari, T. Szamuely: {\it Arithmetic duality theorems 
for 1-motives} , J. Reine Angew. Math.  {\bf 578}, 93--128 (2005).

\bibitem[HSz2]{HS} D. Harari, T. Szamuely: {\it Local-global principles for 1-motives}, Duke Math. J. {\bf 143}, No 3, 531--557 (2008).


\bibitem[Har1]{Harder67} G. Harder: {\it Halbeinfache Gruppenschemata \"uber Dedekindringen}, Invent. Math. {\bf 4}, 165--191 (1967).

\bibitem[Har2]{harder} G. Harder: {\it  \"Uber die Galoiskohomologie halbeinfacher algebraischer Gruppen. III.} J. Reine Angew. Math. {\bf 274--275}, 125--138 (1975).

\bibitem[Izq]{Izq17-SMF} D. Izquierdo: {\it Principe local-global pour les corps de fonctions sur des corps locaux sup\'erieurs II}, Bulletin de la SMF, Vol. {\bf 145} (2), 267--293 (2017).

\bibitem[Mil]{MilADT} J. S. Milne: {\it Arithmetic Duality Theorems}, Second edition,
BookSurge, LLC, Charleston, SC, 2006.

\bibitem[Pra]{prasad} G. Prasad: {\it Strong approximation for
semi-simple groups over function fields}, Ann. of Math. (2) {\bf 105},
no. 3, 553–-572 (1977).

\bibitem[Ray]{Ray} M. Raynaud: {\it Faisceaux amples sur les sch\'emas en groupes et les espaces homog\`enes}, Lecture Notes in Mathematics {\bf 119}, Springer 1970.

\bibitem[Ros1]{RosBig} Z. Rosengarten: {\it Tate duality in positive dimension over function fields}, preprint.

\bibitem[Ros2]{Ros-Tam} Z. Rosengarten: {\it Tamagawa numbers and other invariants of pseudo-reductive groups over global function fields}, to appear in Algebra and Number Theory, 2020.

\bibitem[San]{San} J.-J. Sansuc, {\it Groupe de Brauer et arithm\'etique des groupes alg\'ebriques lin\'eaires sur un corps de nombres},
J. Reine Angew. Math. {\bf 327} (1981) 12–-80.

\bibitem[Tha]{nguyen} N. Q. Thang : {\it On Galois cohomology of semisimple
	groups over local and global fields of positive characteristic, III},
Math. Z. {\bf 275} (2013), no. 3-4, 1287–-1315.

\end{thebibliography}
\end{document}